\documentclass[11pt,english]{article}
\usepackage[osf]{mathpazo}

\usepackage[T1]{fontenc}
\usepackage[utf8]{inputenc}
\usepackage{color}
\usepackage{babel}
\usepackage{amsmath}
\usepackage{amsthm}
\usepackage{amssymb}
\usepackage{stmaryrd}
\usepackage[letterpaper]{geometry}
\usepackage[bookmarks=false,
 breaklinks=false,pdfborder={0 0 1},backref=false,colorlinks=true]
 {hyperref}
\hypersetup{pdftitle={Uniform and Relative Non-Associativity on Minimal Free Resolutions},
 pdfauthor={Michael Nelson},
 linkcolor=blue,citecolor=blue,urlcolor=blue}

\makeatletter

\providecommand{\tabularnewline}{\\}

\theoremstyle{plain}
\newtheorem{thm}{\protect\theoremname}[section]
\theoremstyle{remark}
\newtheorem{rem}[thm]{\protect\remarkname}
\theoremstyle{plain}
\newtheorem{prop}[thm]{\protect\propositionname}
\theoremstyle{remark}
\newtheorem*{rem*}{\protect\remarkname}
\theoremstyle{plain}
\newtheorem{lem}[thm]{\protect\lemmaname}
\theoremstyle{plain}
\newtheorem{cor}[thm]{\protect\corollaryname}
\theoremstyle{definition}
\newtheorem{example}[thm]{\protect\examplename}

\@ifundefined{date}{}{\date{}}
\usepackage{babel}
\usepackage{booktabs}
\usepackage{array}\usepackage{needspace}
\usepackage{tikz}
\usetikzlibrary{calc,arrows.meta}
\usepackage{xurl}
\usepackage{etoolbox}
\AtBeginEnvironment{thebibliography}{\raggedright}
\allowdisplaybreaks[2]

\newtheorem*{mainA}{Theorem A}\newtheorem*{mainB}{Theorem B}

\makeatother

\providecommand{\corollaryname}{Corollary}
\providecommand{\examplename}{Example}
\providecommand{\lemmaname}{Lemma}
\providecommand{\propositionname}{Proposition}
\providecommand{\remarkname}{Remark}
\providecommand{\theoremname}{Theorem}

\begin{document}
\title{Uniform and Relative Non-Associativity\\
on Minimal Free Resolutions}
\author{Michael Nelson}

\maketitle
 
\begin{abstract}
We study unital chain-map multiplications on minimal free resolutions
of quotient rings. We prove that Katthän's example admits no associative
multiplication over any field, without assuming graded commutativity
or compatibility with an internal grading. The proof gives a relative
obstruction in which associativity on a specified finite collection
of triples forces failure at another triple. The argument extends
to a family that includes nonmonomial ideals. For Katthän's ideal,
every Avramov obstruction associated to a regular-sequence quotient
vanishes. We also give criteria for uniform non-associativity, where
a fixed triple fails for every multiplication, by comparing either
associators or their primitives under changes of multiplication. These
criteria persist from monomial elimination ideals and give explicit
obstructions for Avramov's example and the known eight-vertex Bier
sphere over every field. As an application, we show that the edge
ideal of the eight-cycle admits an injective Taylor comparison on
five selected generators, while every selection of that size fails
to admit a multiplicative comparison. Four generators admit a multiplicative
embedding.
\end{abstract}

\section*{Introduction}

Let $(R,\mathfrak{m})$ be a local noetherian ring, and let $F\to R/I$
be a minimal free resolution with $F_{0}=R$. A multiplication on
$F$ is a unital chain map $F\otimes_{R}F\to F$. We call it commutative
when it is strictly graded-commutative, including the vanishing of
odd squares. Buchsbaum and Eisenbud proved that a commutative multiplication
always exists and asked whether an associative one can always be chosen
\cite{BE77}. Counterexamples were subsequently constructed in \cite{Avr81,Sri92,Sri96,Kat19}.

We study how associators change when the multiplication varies. Uniform
non-associativity at a homogeneous triple means that its associator
is nonzero for every multiplication. A relative obstruction imposes
associativity on a prescribed collection of triples and forces its
failure at another triple. Our main result uses this second form of
obstruction to remove both the internal-grading and graded-commutativity
hypotheses from the nonexistence theorem for Katthän's example.

Katthän proved that the minimal resolution of $\Bbbk[x,y,z,w]/\langle x^{2},w^{2},zw,xy,y^{2}z^{2}\rangle$
admits no commutative DG algebra structure whose multiplication respects
the standard internal grading \cite[Theorem 5.1]{Kat19}. We prove
the following statement for all unital chain-map multiplications.

\begin{mainA}[Unrestricted relative non-associativity]Over every
field $\Bbbk$, the minimal resolution over $\Bbbk[x,y,z,w]$ of
\[
\Bbbk[x,y,z,w]/\langle x^{2},w^{2},zw,xy,y^{2}z^{2}\rangle
\]
admits no associative multiplication. Associativity on the $32$ core
triples specified in \eqref{eq:C} forces non-associativity at $(e_{1},e_{5},e_{2})$.
The same assertions hold after localization at $\langle x,y,z,w\rangle$
and completion there.\end{mainA}

Theorems \ref{thm:main} and \ref{thm:relative} prove this statement.
Hancharuk, Laurent-Gengoux, and Strobl compute an explicit $A_{\infty}$-structure
for the same ideal, including the reference associator primitive with
coefficient $yz$ \cite[Example 4.3]{HLS24}. Our contribution is
the obstruction to simultaneous cancellation of the associators under
arbitrary changes of multiplication. The familiar reference associator
at $(e_{1},e_{5},e_{2})$ can be cancelled by a change of multiplication,
as Example \ref{ex:cancel} shows. Associativity on the core triples
constrains the same correction coefficients. Reduction modulo $\langle x,w\rangle$
and Nakayama's lemma then force a nonzero associator at the distinguished
triple. Theorem \ref{thm:family} extends this argument by replacing
$y^{2}z^{2}$ with $yz\,c(y,z)$ for any nonzero $c\in\langle y,z\rangle$,
including nonhomogeneous polynomials.

Katthän also noted that he could not find a regular sequence in this
ideal to which Avramov's obstruction applies \cite[Remark~5.4]{Kat19}.
Theorem~\ref{thm:katthan-avramov} explains this difficulty: every
Avramov obstruction associated to a regular-sequence quotient vanishes,
over every field. The result includes nonhomogeneous regular sequences
and remains valid after localization and completion. Thus the relative
obstruction detects a failure of associativity that these change-of-rings
obstructions do not detect. 

The uniform criteria compare a reference associator with all of its
homotopy corrections. The first criterion detects a class modulo $\mathfrak{m}^{2}F$
and endpoint products with boundaries.

\begin{mainB}[Uniform non-associativity]Let $(R,\mathfrak{m})$ be
local and noetherian, and let $F$ be a minimal free resolution of
$R/I$ with reference multiplication $\mu$. If homogeneous $a_{1},a_{2},a_{3}\in F$
satisfy
\[
\mathrm{d}a_{i}\in\mathfrak{m}^{2}F,\qquad a_{1}a_{2},a_{2}a_{3}\in\mathfrak{m}F,\qquad[a_{1},a_{2},a_{3}]_{\mu}\notin\mathfrak{m}^{2}F+a_{1}\mathrm{d}F+a_{3}\mathrm{d}F,
\]
then every multiplication on $F$ is non-associative at this triple.
For monomial ideals, these hypotheses persist from an elimination
ideal to the original ideal when the reference multiplication is multigraded.\end{mainB}

This combines Theorem \ref{theoremnonassociativeforall} with Proposition
\ref{prop:elimination-persistence}. For degree-one inputs, Theorem
\ref{thm:primitive-criterion} gives a stronger criterion by testing
a primitive of the associator one homological degree higher. The primitive
is compared modulo coefficient ideals, boundaries, and endpoint products.
Proposition \ref{prop:error-ideals} gives coefficient-ideal estimates
in arbitrary homological degrees. Powers of Avramov's example yield
uniform obstructions with arbitrarily high coefficient order.

The elimination argument applies these criteria to the known eight-vertex
Bier sphere, identifying an explicit fixed triple at which every multiplication
fails over every field. Denham and Suciu studied this sphere and its
Massey-product obstruction in the characteristic-zero setting of \cite[Section 8]{DS07},
within Baskakov's family \cite{Bas03}. Earlier grade-five Gorenstein
counterexamples appear in \cite{Sri96}. The permanence argument uses
the retained-resolution construction of \cite{Boo12} and is related
to DG-sensitive pruning \cite{GMP26}.

Uniformity also obstructs multiplicative Taylor comparisons when both
the target multiplication and the comparison map vary. For the edge
ideal of the eight-cycle, an injective Taylor comparison on five selected
generators exists, but every selection of that size fails to admit
a multiplicative comparison. Four generators admit a multiplicative
embedding. Section \ref{sec:tcomp} proves this sharp result and a
family of related cycle examples.

Section \ref{sec:uniform} develops the criteria and their first examples.
Section \ref{sec:examples} proves permanence under elimination and
gives the path and Gorenstein applications. Section \ref{sec:katthan}
proves the relative obstruction, its family extension, and the vanishing
of the regular-sequence Avramov obstructions, and Section \ref{sec:tcomp}
treats Taylor comparisons.

Throughout, degrees of maps and graded commutativity refer to the
homological grading. For polynomial rings, $\mathfrak{m}$ denotes
the homogeneous maximal ideal. Competing multiplications and homotopies
may have arbitrary internal degrees unless stated otherwise.

\section{Criteria for uniform non-associativity}

\label{sec:uniform}

Throughout this section, let $R$ be a commutative ring, let $I\subsetneq R$
be an ideal, and let $F\to R/I$ be a nonnegative free resolution
over $R$ with $F_{0}=R$.

\subsection{Multiplications and generalized associators}

Recall that by a multiplication on $F$ we mean a unital chain map
$\mu\colon F\otimes_{R}F\to F$, not necessarily associative. We call
it a commutative multiplication when it is also strictly graded-commutative.
Every multiplication on $F$ lifts the usual multiplication on $R/I$:
$R$-bilinearity and unitality give $\mu(r,s)=rs\mu(1,1)=rs$ for
$r,s\in R$, and in positive degrees the augmentation is zero. For
an arbitrary chain map $\alpha\colon F\otimes_{R}F\to F$, the lifting
condition is instead $\alpha(1,1)-1\in I$; it is not automatic.

In particular, after fixing a multiplication $\mu$, every other multiplication
on $F$ has the form $\mu_{h}:=\mu+\mathrm{d}h+h\mathrm{d}$, where
$h\colon F\otimes_{R}F\to F$ is a homotopy, that is, a graded $R$-linear
map of homological degree one. This is the usual comparison-homotopy
argument: $F\otimes_{R}F$ is a bounded-below complex of free modules,
and the two maps have the same composite with the augmentation $F\to R/I$.
No normalization of $h$ is needed for the first criterion below.

For minimal free resolutions over local noetherian rings, Buchsbaum
and Eisenbud showed in \cite{BE77} that a commutative multiplication
always exists. They asked whether an associative commutative multiplication
can be chosen. It was subsequently shown that this need not be the
case; see \cite{Avr81,Kat19,Sri92,Sri96} for counterexamples. In
this section, we show that if a multiplication $\mu$ fails to be
associative at a homogeneous triple $(a_{1},a_{2},a_{3})$ and some
additional conditions hold for $\mu$ and that triple, then every
multiplication on $F$ fails to be associative at the same triple.

Before we state and prove the theorem, we set up notation. Equip $F$
with a fixed multiplication $\mu$, and write $a_{1}a_{2}=\mu(a_{1},a_{2})$.
The associator for $\mu$ is the chain map $[\cdot]_{\mu}\colon F\otimes_{R}F\otimes_{R}F\to F$
given by $[\cdot]_{\mu}=\mu(\mu\otimes1-1\otimes\mu)$. The corresponding
$R$-trilinear map is denoted $[\cdot,\cdot,\cdot]_{\mu}\colon F^{3}\to F$.
In particular, for $a_{1},a_{2},a_{3}\in F$,
\[
[a_{1}\otimes a_{2}\otimes a_{3}]_{\mu}=[a_{1},a_{2},a_{3}]_{\mu}=(a_{1}a_{2})a_{3}-a_{1}(a_{2}a_{3}).
\]

Since $\mu$ is fixed, we simplify notation by dropping its subscript
when there is no ambiguity. If $\mu_{h}=\mu+\mathrm{d}h+h\mathrm{d}$
is another multiplication, write $a\star b=\mu_{h}(a,b)$. Its associator
is given by
\begin{equation}
[\cdot]_{\mu_{h}}=[\cdot]+\mathrm{d}H+H\mathrm{d},\label{eq:homotopy-associator}
\end{equation}
where $H=[\cdot]_{\mu,h}+[\cdot]_{h,\mu_{h}}$ and
\begin{equation}
[\cdot]_{\mu,h}=\mu(h\otimes1-1\otimes h),\qquad[\cdot]_{h,\mu_{h}}=h(\mu_{h}\otimes1-1\otimes\mu_{h}).\label{eq:generalized-associators}
\end{equation}

We refer to these maps as generalized associators because they generalize
the usual associator in the sense that $[\cdot]_{\mu,\mu}=[\cdot]_{\mu}$.
Additional signs appear when evaluating \eqref{eq:generalized-associators}
on elements, according to the Koszul sign rule. In \eqref{eq:homotopy-associator},
we use the same symbol $\mathrm{d}$ for the differentials on $F$
and $F^{\otimes3}$, with the domain making its meaning clear; for
instance, the differential in $\mathrm{d}H$ is that of $F$ because
$H$ lands in $F$. We use the same convention for the differential
on $F^{\otimes2}$.

Finally, we can decompose the generalized associator $[\cdot]_{h,\mu_{h}}$
into three generalized associators:
\[
[\cdot]_{h,\mu_{h}}=[\cdot]_{h,\mu}+[\cdot]_{h,\mathrm{d}h}+[\cdot]_{h,h\mathrm{d}}.
\]
The notation indicates how these maps are defined; for example, $[\cdot]_{h,h\mathrm{d}}=h(h\mathrm{d}\otimes1-1\otimes h\mathrm{d})$.
With this notation in place, we are ready to state and prove the criterion.

\Needspace{14\baselineskip}

\subsection{The uniform criterion}
\begin{thm}[Uniform non-associativity]
\label{theoremnonassociativeforall} Let $(R,\mathfrak{m})$ be a
local noetherian ring and let $F$ be a minimal free resolution of
$R\slash I$ over $R$, equipped with a multiplication $\mu$. Suppose
there are homogeneous elements $a_{1},a_{2},a_{3}\in F$ such that
\[
\mathrm{d}a_{1},\mathrm{d}a_{2},\mathrm{d}a_{3}\in\mathfrak{m}^{2}F,\qquad a_{1}a_{2},a_{2}a_{3}\in\mathfrak{m}F,
\]
and 
\[
[a_{1},a_{2},a_{3}]\notin\mathfrak{m}^{2}F+a_{1}\mathrm{d}F+a_{3}\mathrm{d}F.
\]
Then every multiplication on $F$ fails to be associative at $(a_{1},a_{2},a_{3})$. 
\end{thm}

\begin{proof}
Let $\mu_{h}=\mu+\mathrm{d}h+h\mathrm{d}$ be any competing multiplication,
by the comparison-homotopy argument above. We analyze the correction
$\mathrm{d}H+H\mathrm{d}$ in \eqref{eq:homotopy-associator} modulo
$\mathfrak{m}^{2}F$. Set $a_{1,2}=h(a_{1},a_{2})$ and $a_{2,3}=h(a_{2},a_{3})$.
We have 
\begin{align*}
(\mathrm{d}H+H\mathrm{d})(a_{1}\otimes a_{2}\otimes a_{3}) & \equiv\mathrm{d}H(a_{1}\otimes a_{2}\otimes a_{3})\\
 & =\mathrm{d}\bigl([\cdot]_{\mu,h}+[\cdot]_{h,\mu}+[\cdot]_{h,\mathrm{d}h}+[\cdot]_{h,h\mathrm{d}}\bigr)(a_{1}\otimes a_{2}\otimes a_{3})\\
 & \equiv\mathrm{d}\bigl([\cdot]_{\mu,h}+[\cdot]_{h,\mu}\bigr)(a_{1}\otimes a_{2}\otimes a_{3})\\
 & \equiv\mathrm{d}[\cdot]_{\mu,h}(a_{1}\otimes a_{2}\otimes a_{3})\\
 & =\mathrm{d}\bigl(a_{1,2}a_{3}-(-1)^{|a_{1}|}a_{1}a_{2,3}\bigr)\\
 & \equiv(\mathrm{d}a_{1,2})a_{3}-a_{1}\mathrm{d}a_{2,3}.
\end{align*}
The first congruence uses $\mathrm{d}a_{i}\in\mathfrak{m}^{2}F$.
For the second, minimality gives $\mathrm{d}F\subseteq\mathfrak{m}F$,
and the differential occurs twice in each of the discarded terms.
The third uses $a_{1}a_{2},a_{2}a_{3}\in\mathfrak{m}F$, followed
by the outer differential. The last congruence follows from the Leibniz
rule and $\mathrm{d}a_{1},\mathrm{d}a_{3}\in\mathfrak{m}^{2}F$.

Consequently, 
\begin{equation}
[a_{1},a_{2},a_{3}]_{\mu_{h}}\equiv[a_{1},a_{2},a_{3}]+(\mathrm{d}a_{1,2})a_{3}-a_{1}\mathrm{d}a_{2,3}\pmod{\mathfrak{m}^{2}F}.\label{eq:uniform-congruence}
\end{equation}
Let $\mu^{\mathrm{op}}$ be the graded opposite multiplication, defined
by 
\[
\mu^{\mathrm{op}}(a,b)=(-1)^{|a||b|}\mu(b,a).
\]
The graded twist is a chain map, so $\mu^{\mathrm{op}}$ is again
a multiplication. The comparison-homotopy argument and minimality
give the first inclusion below; $R$-bilinearity and $\mathrm{d}F\subseteq\mathfrak{m}F$
give the second: 
\begin{equation}
(\mu-\mu^{\mathrm{op}})(F,F)\subseteq\mathfrak{m}F,\qquad(\mu-\mu^{\mathrm{op}})(\mathrm{d}F,F)\subseteq\mathfrak{m}^{2}F.\label{eq:graded-opposite-order}
\end{equation}
Thus the correction in \eqref{eq:uniform-congruence} lies in $\mathfrak{m}^{2}F+a_{1}\mathrm{d}F+a_{3}\mathrm{d}F$.
The nonmembership hypothesis therefore gives $[a_{1},a_{2},a_{3}]_{\mu_{h}}\ne0$. 
\end{proof}
\begin{rem}
\label{rem:scope-uniform} Neither the reference multiplication nor
the competing multiplication is required to be commutative. 

For a minimal graded resolution over a polynomial ring, the criterion
can be applied after localization at the homogeneous maximal ideal.
Every polynomial multiplication, whether internally graded or not,
localizes to a multiplication there. 
\end{rem}

\subsection{Avramov's example}

\label{sec:avramov-uniform}

We apply Theorem~\ref{theoremnonassociativeforall} to Avramov's
example \cite[Example~2.2]{Avr81}. Avramov proves that its minimal
resolution admits no DG module structure over the Koszul complex on
the two squared generators. The Massey-product interpretation of this
obstruction appears in \cite[Example~5.1.4]{Avr81}, and the resulting
nonexistence of a commutative DG algebra structure is stated explicitly
in \cite[Remark~5.2.4]{Avr81}. We will show that one fixed triple
is non-associative for every multiplication on the resolution.

Let $R=\Bbbk[x,y,z,w]$ be a polynomial ring over a field $\Bbbk$,
let $\mathfrak{m}=\langle x,y,z,w\rangle$, and put 
\[
I=\langle x^{2},w^{2},zw,xy,yz\rangle.
\]
Write $F$ for the minimal free resolution of $R/I$ over $R$. We
use the conventions for labeled cellular resolutions of \cite{BS98};
see also \cite[Chapter~4]{MS05}. One can visualize $F$ as being
supported on the $I$-labeled cellular complex below: 
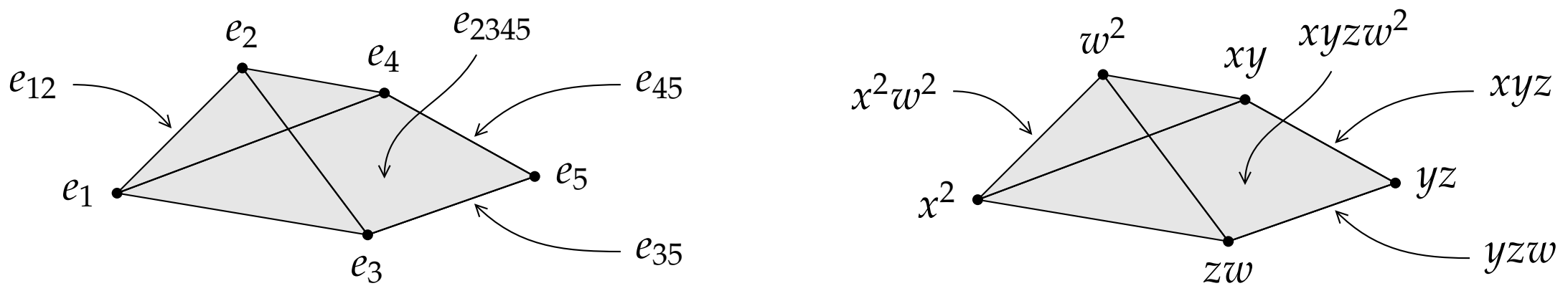
\begin{center}
\begin{tabular}{cccc}
\begin{tikzpicture}[scale=0.77]

\draw[fill=gray!20] (0,0) -- (3,-0.5) -- (5,0.2) -- (3.2,1.2)-- (0,0); 
\draw[fill=gray!20] (0,0) -- (1.5,1.5) -- (3.2,1.2)-- (0,0); 
\draw[] (1.5,1.5) -- (3,-0.5); \draw[fill=gray!20] (3,-0.5) -- (5,0.2) -- (3.2,1.2);
\draw[color=black!100] (1.5,1.5) -- (3,-0.5) -- (5,0.2);
\draw[color=black!100] (0,0) -- (3.2,1.2) -- (5,0.2);

\node[circle, fill=black, inner sep=1pt, label=left:$e_1 $] (a) at (0,0) {};
\node[circle, fill=black, inner sep=1pt, label=above:$e_2 $] (b) at (1.5,1.5) {};
\node[circle, fill=black, inner sep=1pt, label=below:$e_3 $] (c) at (3,-0.5) {};
\node[circle, fill=black, inner sep=1pt, label=above:$e_4 $] (d) at (3.2,1.2) {}; 
\node[circle, fill=black, inner sep=1pt, label=right:$e_5 $] (e) at (5,0.2) {};

\node[] (x) at (-1,1.3) {$ e_{12} $}; \node[] (x') at (0.8,0.6) {};
\draw[-{Straight Barb[length=3pt,width=3pt]}] (x) edge[out=0, in=130] node[below] {} (x');

\node[] (y) at (6.5,1.3) {$ e_{45} $}; \node[] (y') at (4.1,0.5) {}; 
\draw[-{Straight Barb[length=3pt,width=3pt]}] (y) edge[out=180, in=45] node[below] {} (y');

\node[] (z) at (6.5,-0.7) {$ e_{35} $}; \node[] (z') at (4.1,0.05) {};
\draw[-{Straight Barb[length=3pt,width=3pt]}] (z) edge[out=180, in=-45] node[below] {} (z');

\node[] (w) at (4.5,2) {$ e_{2345} $}; \node[] (w') at (3.2,0) {}; 
\draw[-{Straight Barb[length=3pt,width=3pt]}] (w) edge[out=-120, in=90] node[below] {} (w');

\end{tikzpicture} &  &  & \begin{tikzpicture}[scale=0.77]

\draw[fill=gray!20] (0,0) -- (3,-0.5) -- (5,0.2) -- (3.2,1.2)-- (0,0); 
\draw[fill=gray!20] (0,0) -- (1.5,1.5) -- (3.2,1.2)-- (0,0); 
\draw[] (1.5,1.5) -- (3,-0.5); \draw[fill=gray!20] (3,-0.5) -- (5,0.2) -- (3.2,1.2);
\draw[color=black!100] (1.5,1.5) -- (3,-0.5) -- (5,0.2);
\draw[color=black!100] (0,0) -- (3.2,1.2) -- (5,0.2);

\node[circle, fill=black, inner sep=1pt, label=left:$x^2 $] (a) at (0,0) {};
\node[circle, fill=black, inner sep=1pt, label=above:$w^2 $] (b) at (1.5,1.5) {};
\node[circle, fill=black, inner sep=1pt, label=below:$zw $] (c) at (3,-0.5) {};
\node[circle, fill=black, inner sep=1pt, label=above:$xy $] (d) at (3.2,1.2) {}; 
\node[circle, fill=black, inner sep=1pt, label=right:$yz $] (e) at (5,0.2) {};

\node[] (x) at (-1,1.3) {$ x^2 w^2 $}; \node[] (x') at (0.8,0.6) {};
\draw[-{Straight Barb[length=3pt,width=3pt]}] (x) edge[out=0, in=130] node[below] {} (x');

\node[] (y) at (6.5,1.3) {$ x y z $}; \node[] (y') at (4.1,0.5) {}; 
\draw[-{Straight Barb[length=3pt,width=3pt]}] (y) edge[out=180, in=45] node[below] {} (y');

\node[] (z) at (6.5,-0.7) {$ y z w $}; \node[] (z') at (4.1,0.05) {};
\draw[-{Straight Barb[length=3pt,width=3pt]}] (z) edge[out=180, in=-45] node[below] {} (z');

\node[] (w) at (4.5,2) {$ x y z w^2 $}; \node[] (w') at (3.2,0) {}; 
\draw[-{Straight Barb[length=3pt,width=3pt]}] (w) edge[out=-120, in=90] node[below] {} (w');

\end{tikzpicture}\tabularnewline
\end{tabular}
\par\end{center}

We write down the homogeneous components of $F$ as a graded $R$-module:
\begin{align*}
F_{0} & =R,\\
F_{1} & =Re_{1}+Re_{2}+Re_{3}+Re_{4}+Re_{5},\\
F_{2} & =Re_{12}+Re_{13}+Re_{14}+Re_{23}+Re_{24}+Re_{35}+Re_{45},\\
F_{3} & =Re_{123}+Re_{124}+Re_{1345}+Re_{2345},\\
F_{4} & =Re_{12345}.
\end{align*}
The vertex labels, in order, are $x^{2},w^{2},zw,xy,yz$. We use compatible
cellular orientations: on the simplicial faces the differential is
the usual signed least-common-multiple differential, and, writing
$\tau=e_{12345}$, we have 
\begin{align}
\mathrm{d}e_{1345} & =x^{2}e_{35}-xwe_{45}-zwe_{14}+ye_{13},\label{eq:av-cell-1345}\\
\mathrm{d}e_{2345} & =xwe_{35}-w^{2}e_{45}-ze_{24}+xye_{23},\label{eq:av-cell-2345}\\
\mathrm{d}\tau & =-ye_{123}+ze_{124}-we_{1345}+xe_{2345}.\label{eq:av-top-differential}
\end{align}
For example, $\mathrm{d}e_{123}=x^{2}e_{23}-we_{13}+ze_{12}$ and
$\mathrm{d}e_{124}=xe_{24}-w^{2}e_{14}+ye_{12}$.

Choose a multigraded commutative multiplication $\mu$ on this resolution;
such a choice is obtained by carrying out the usual lifting construction
multidegree by multidegree. Multigrading and the Leibniz rule force
the products needed below: 
\begin{align}
e_{1}e_{5} & =ze_{14}+xe_{45}, & e_{2}e_{5} & =ye_{23}+we_{35},\label{eq:av-basic-products}\\
e_{2}e_{14} & =-e_{124}, & e_{1}e_{23} & =e_{123},\nonumber \\
e_{2}e_{45} & =-e_{2345}, & e_{1}e_{35} & =e_{1345}.\nonumber 
\end{align}
These formulas give 
\begin{equation}
[e_{1},e_{5},e_{2}]_{\mu}=ye_{123}-ze_{124}+we_{1345}-xe_{2345}=-\mathrm{d}\tau.\label{eq:av-linear-associator}
\end{equation}
The coefficients of this associator are linear. We verify that it
satisfies all the hypotheses of Theorem~\ref{theoremnonassociativeforall}.
\begin{prop}
\label{prop:avramov-uniform} Every multiplication $F\otimes_{R}F\to F$
is non-associative at $(e_{1},e_{5},e_{2})$. In particular, the conclusion
holds for every commutative multiplication, without any assumption
on its internal grading. 
\end{prop}

\begin{proof}
All three vertex differentials $x^{2},yz,w^{2}$ belong to $\mathfrak{m}^{2}$.
By \eqref{eq:av-basic-products} and graded commutativity, both consecutive
products $e_{1}e_{5}$ and $e_{5}e_{2}$ belong to $\mathfrak{m}F_{2}$.
It remains to check the boundary-product terms in the criterion. The
endpoint products with $F_{3}$ are 
\[
\begin{array}{c|rrrr}
 & e_{123} & e_{124} & e_{1345} & e_{2345}\\
\hline e_{1}\cdot & 0 & 0 & 0 & x\tau\\
e_{2}\cdot & 0 & 0 & -w\tau & 0.
\end{array}
\]
Thus $e_{1}F_{3},e_{2}F_{3}\subseteq\mathfrak{m}F_{4}$. For $i\in\{1,2\}$
and $C\in F_{3}$, Leibniz gives 
\[
e_{i}\,\mathrm{d}C=(\mathrm{d}e_{i})C-\mathrm{d}(e_{i}C)\in\mathfrak{m}^{2}F_{3},
\]
because $\mathrm{d}e_{i}\in\mathfrak{m}^{2}$ and the resolution is
minimal. Consequently 
\[
\bigl(\mathfrak{m}^{2}F+e_{1}\mathrm{d}F+e_{2}\mathrm{d}F\bigr)_{3}=\mathfrak{m}^{2}F_{3}.
\]
The element $-\mathrm{d}\tau$ in \eqref{eq:av-linear-associator}
has nonzero linear coefficients, so it does not belong to this submodule.
Apply Theorem~\ref{theoremnonassociativeforall} after localizing
at $\mathfrak{m}$; the nonzero linear class remains nonzero, and
any polynomial multiplication localizes.

Remark~\ref{rem:scope-uniform} shows that the competing multiplication
need not be graded-commutative. In fact, the calculation above and
\eqref{eq:uniform-congruence} give the stronger congruence 
\[
[e_{1},e_{5},e_{2}]_{\nu}\equiv-\mathrm{d}\tau\not\equiv0\pmod{\mathfrak{m}^{2}F_{3}}.
\]
\end{proof}

\subsection{A criterion using a primitive of the associator}

\label{sec:further-criteria}

Fix a reference multiplication $\mu\colon F\otimes_{R}F\to F$, write
$ab=\mu(a,b)$, and write $a\star b=\mu_{h}(a,b)$ for a competing
multiplication.

Theorem \ref{theoremnonassociativeforall} detects an associator modulo
$\mathfrak{m}^{2}F$ and endpoint products with boundaries. For degree-one
triples, we can instead test a primitive one homological degree higher.
Normalizing the comparison homotopy makes the terms involving the
differentials of the three inputs vanish exactly, so no hypothesis
$\mathrm{d}a_{i}\in\mathfrak{m}^{2}F$ is needed. The condition $\mathrm{d}F_{3}\subseteq JF_{2}$
will control the boundaries of $h(a_{1},a_{2})$ and $h(a_{2},a_{3})$;
for a minimal resolution and $J=\mathfrak{m}$, this condition is
automatic. Both criteria establish non-associativity at one fixed
triple for every multiplication.

Let $a_{1},a_{2},a_{3}\in F$ be homogeneous, put $s=|a_{1}|+|a_{2}|+|a_{3}|>0$,
and suppose that $[a_{1},a_{2},a_{3}]$ is a cycle. Exactness gives
$\tau\in F_{s+1}$ such that 
\begin{equation}
\mathrm{d}\tau=[a_{1},a_{2},a_{3}].\label{eq:associator-primitive}
\end{equation}
We call $\tau$ a \emph{primitive of the associator}. Two choices
differ by an element of $\mathrm{d}F_{s+2}$.

If $a_{1},a_{2},a_{3}\in F_{1}$, the associator is automatically
a cycle. Indeed, unitality and $R$-bilinearity put $R=F_{0}$ in
the nucleus, meaning that every associator with a scalar input vanishes.
The Leibniz rule for the associator then gives 
\[
\mathrm{d}[a_{1},a_{2},a_{3}]=[r_{1},a_{2},a_{3}]-[a_{1},r_{2},a_{3}]+[a_{1},a_{2},r_{3}]=0,
\]
where we set $r_{i}=\mathrm{d}a_{i}$.
\begin{rem*}
The primitive formulation has a natural interpretation in terms of
$A_{\infty}$-structures. We use homological conventions: $\mu_{n}\colon F^{\otimes n}\to F$
has degree $n-2$, and the differential on the Hom complex is 
\[
\mathrm{d}^{\star}(\varphi)=\mathrm{d}\varphi-(-1)^{|\varphi|}\varphi\mathrm{d}.
\]
With $\mu_{2}=\mu$, we choose the sign convention for $\mu_{3}$
so that the arity-three identity is 
\[
\mathrm{d}^{\star}(\mu_{3})=\mathrm{d}\mu_{3}+\mu_{3}\mathrm{d}=[\cdot]_{\mu}.
\]
One can choose $\mu_{3}$ to vanish whenever an input is $1$, by
the normalized lifting argument of Lemma~\ref{lem:normalized}. For
$a_{1},a_{2},a_{3}\in F_{1}$, this gives 
\[
\mathrm{d}\mu_{3}(a_{1},a_{2},a_{3})=[a_{1},a_{2},a_{3}]_{\mu}.
\]
Thus $\mu_{3}(a_{1},a_{2},a_{3})$ is a choice of the primitive in
\eqref{eq:associator-primitive}. The existence and transfer of strictly
unital $A_{\infty}$-structures on free resolutions over $R$ are
standard; see \cite[Section~5, especially Proposition~5.1]{BCLP25}.
\end{rem*}
Recall that for a reference chain map $\mu\colon F\otimes_{R}F\to F$
lifting the multiplication on $R/I$, every other such chain map has
the form $\mu_{h}=\mu+\mathrm{d}h+h\mathrm{d}$ for some homotopy
$h$. The lifting condition alone does not guarantee that $h$ can
be normalized. For multiplications, the two maps agree on scalar inputs,
and the homotopy can be chosen to vanish there as well.
\begin{lem}[Normalized comparison]
\label{lem:normalized}Let $\mu$ be a multiplication on $F$. Every
other multiplication can be written as $\mu_{h}=\mu+\mathrm{d}h+h\mathrm{d}$,
where $h\colon F\otimes_{R}F\to F$ is a graded $R$-linear map of
homological degree one satisfying
\[
h(r,a)=h(a,r)=0\qquad(r\in R,\ a\in F).
\]
If both $\mu$ and $\mu_{h}$ are commutative multiplications, then
$h$ can additionally be chosen to satisfy
\[
h(a,b)=(-1)^{|a||b|}h(b,a),\qquad h(a,a)=0\quad\text{if }|a|\text{ is odd}
\]
for homogeneous $a,b\in F$.
\end{lem}

\begin{proof}
Put $\Delta=\mu_{h}-\mu$ and $U=F_{0}\otimes_{R}F+F\otimes_{R}F_{0}$.
The submodule $U$ is a subcomplex of $F\otimes_{R}F$, and unitality
implies that $\Delta$ vanishes on $U$. Thus $\Delta$ factors through
a chain map $\overline{\Delta}\colon C=(F\otimes_{R}F)/U\to F$. The
module $C_{n}$ is the direct sum of $F_{p}\otimes_{R}F_{q}$ over
$p,q>0$ with $p+q=n$; in particular it is free, and $C_{0}=C_{1}=0$.

Let $\epsilon\colon F\to R/I$ be the augmentation. Since $C_{0}=0$,
we have $\epsilon\overline{\Delta}=0$. The bounded-below complex
$C$ is degreewise free, so the homotopy lifting argument gives a
graded map $\overline{h}\colon C\to F$ of degree one with
\[
\overline{\Delta}=\mathrm{d}\overline{h}+\overline{h}\mathrm{d}_{C}.
\]
Composing $\overline{h}$ with the quotient map $F\otimes_{R}F\to C$
gives the required normalized homotopy.

If both multiplications are commutative, let $V\subseteq F\otimes_{R}F$
be the graded submodule generated by
\[
a\otimes b-(-1)^{|a||b|}b\otimes a,\qquad a\otimes a\quad\text{if }|a|\text{ is odd},
\]
where $a,b$ are homogeneous. It is a subcomplex: the symmetry relations
are preserved by the differential, and for odd $a$, $\mathrm{d}(a\otimes a)=\mathrm{d}a\otimes a-a\otimes\mathrm{d}a$
is a symmetry relation. Now $\Delta$ factors through $D=(F\otimes_{R}F)/(U+V)$.
A homogeneous basis of $F$ containing $1$ gives a basis of $D$
indexed by unordered pairs of positive-degree basis elements, omitting
repeated odd basis elements. Thus $D$ is degreewise free and $D_{0}=D_{1}=0$.
The same homotopy lifting argument gives a homotopy vanishing on both
$U$ and $V$, proving the additional identities without division
by two.
\end{proof}
\begin{thm}[Uniform non-associativity from a primitive]
\label{thm:primitive-criterion}Let $a_{1},a_{2},a_{3}\in F_{1}$,
and let $J\subseteq R$ be an ideal such that
\[
a_{1}a_{2},a_{2}a_{3}\in JF_{2},\qquad\mathrm{d}F_{3}\subseteq JF_{2}.
\]
Set 
\[
N=JF_{4}+\mathrm{d}F_{5}+F_{3}a_{3}+a_{1}F_{3}\subseteq F_{4}.
\]
For every multiplication $\mu_{h}$ and every primitive $\tau_{h}$
of its associator at $(a_{1},a_{2},a_{3})$, one has
\[
\tau_{h}\equiv\tau\pmod N.
\]
In particular, if $\tau\notin N$, then every multiplication on $F$
is non-associative at $(a_{1},a_{2},a_{3})$. No hypothesis on the
orders of $\mathrm{d}a_{i}$ is required.
\end{thm}

\begin{proof}
Let $\mu_{h}$ be a competing multiplication, and choose a normalized
homotopy $h$ with $\mu_{h}=\mu+\mathrm{d}h+h\mathrm{d}$ by Lemma
\ref{lem:normalized}. Use the associator homotopy $H=[\cdot]_{\mu,h}+[\cdot]_{h,\mu_{h}}$
from \eqref{eq:homotopy-associator}. Since $\mathrm{d}a_{i}\in F_{0}=R$,
normalization and unitality give

\[
h\mathrm{d}(a_{i}\otimes a_{j})=0,\qquad H\mathrm{d}(a_{1}\otimes a_{2}\otimes a_{3})=0.
\]
Thus $\tau_{h}-\tau-H(a_{1},a_{2},a_{3})$ is a cycle in $F_{4}$,
hence belongs to $\mathrm{d}F_{5}\subseteq N$. It remains to show
that $H(a_{1},a_{2},a_{3})\in N$.

Put $u=h(a_{1},a_{2})$ and $v=h(a_{2},a_{3})$, both in $F_{3}$.
The pairwise vanishing above gives
\[
a_{1}\star a_{2}=a_{1}a_{2}+\mathrm{d}u\in JF_{2},\qquad a_{2}\star a_{3}=a_{2}a_{3}+\mathrm{d}v\in JF_{2}.
\]
In generalized-associator notation, $[a_{1},a_{2},a_{3}]_{h,h\mathrm{d}}=0$,
so
\begin{equation}
\begin{split}H(a_{1},a_{2},a_{3})={} & [a_{1},a_{2},a_{3}]_{\mu,h}\\
 & +[a_{1},a_{2},a_{3}]_{h,\mu}+[a_{1},a_{2},a_{3}]_{h,\mathrm{d}h}.
\end{split}
\label{eq:primitive-variation}
\end{equation}
Each term belongs to the corresponding summand of $N$:
\[
\begin{aligned}{}[a_{1},a_{2},a_{3}]_{\mu,h} & =ua_{3}+a_{1}v\in F_{3}a_{3}+a_{1}F_{3},\\{}
[a_{1},a_{2},a_{3}]_{h,\mu} & =h(a_{1}a_{2},a_{3})-h(a_{1},a_{2}a_{3})\in JF_{4},\\{}
[a_{1},a_{2},a_{3}]_{h,\mathrm{d}h} & =h(\mathrm{d}u,a_{3})-h(a_{1},\mathrm{d}v)\in JF_{4}.
\end{aligned}
\]
The last two inclusions use the hypotheses and $R$-linearity of $h$,
with $h(F_{2}\otimes_{R}F_{1}),h(F_{1}\otimes_{R}F_{2})\subseteq F_{4}$.
This proves the congruence. If $\mu_{h}$ were associative at the
triple, we could choose $\tau_{h}=0$, contradicting $\tau\notin N$.
\end{proof}
\begin{cor}[The residue-field criterion]
\label{cor:primitive-residue} Suppose $(R,\mathfrak{m},\Bbbk)$
is local and $F$ is minimal. Put $V=F/\mathfrak{m}F$, and use a
bar for reduction modulo $\mathfrak{m}$. If $a_{1},a_{2},a_{3}\in F_{1}$
satisfy 
\[
a_{1}a_{2},a_{2}a_{3}\in\mathfrak{m}F_{2},\qquad\bar{\tau}\ne0\quad\text{in}\quad\frac{V_{4}}{V_{3}\bar{a}_{3}+\bar{a}_{1}V_{3}},
\]
then every multiplication on $F$ is non-associative at $(a_{1},a_{2},a_{3})$. 
\end{cor}

\begin{proof}
Apply Theorem~\ref{thm:primitive-criterion} with $J=\mathfrak{m}$.
Minimality gives $\mathrm{d}F_{3}\subseteq\mathfrak{m}F_{2}$ and
$\mathrm{d}F_{5}\subseteq\mathfrak{m}F_{4}$. 
\end{proof}
The product on $V$ is independent of the chosen multiplication on
$F$, since $\mathrm{d}h+h\mathrm{d}$ has image in $\mathfrak{m}F$.
It is the usual product on $\operatorname{Tor}^{R}(R/I,\Bbbk)$. In
particular, the quotient in Corollary~\ref{cor:primitive-residue}
is independent of the reference map. In the general ideal formulation
of Theorem \ref{thm:primitive-criterion}, the submodule $N$ is defined
using the fixed reference multiplication $ab$. Its endpoint products
are taken with respect to this multiplication throughout the comparison.

The endpoint-product quotient in Corollary \ref{cor:primitive-residue}
has the form of the familiar indeterminacy of a triple Massey product
\cite[Section 5.1]{DS07}. The primitive theorem formulates the comparison
over $R$, keeping track of coefficient ideals and endpoint submodules
while fixing the triple and allowing every competing multiplication.
\begin{rem}
\label{rem:primitive-recovers-uniform} For degree-one triples, Corollary~\ref{cor:primitive-residue}
recovers Theorem~\ref{theoremnonassociativeforall}. Indeed, suppose
that the hypotheses of that theorem hold. If 
\[
\tau\in\mathfrak{m}F_{4}+F_{3}a_{3}+a_{1}F_{3},
\]
then applying $\mathrm{d}$, using minimality and $\mathrm{d}a_{1},\mathrm{d}a_{3}\in\mathfrak{m}^{2}F_{0}$,
gives 
\[
[a_{1},a_{2},a_{3}]_{\mu}\in\mathfrak{m}^{2}F_{3}+(\mathrm{d}F_{3})a_{3}+a_{1}\mathrm{d}F_{3}\subseteq\mathfrak{m}^{2}F+a_{1}\mathrm{d}F+a_{3}\mathrm{d}F.
\]
The last inclusion follows from \eqref{eq:graded-opposite-order},
so this contradicts the nonmembership hypothesis. This also shows
that the hypothesis $\mathrm{d}a_{2}\in\mathfrak{m}^{2}F_{0}$ can
be omitted from Theorem~\ref{theoremnonassociativeforall} when all
three inputs have degree one. The original theorem continues to apply
to triples of arbitrary homogeneous degrees. The converse fails even
for degree-one triples: in Example \ref{ex:powered-avramov}, the
displayed associator belongs to $\mathfrak{m}^{2}F_{3}$, while its
primitive has nonzero residue in the quotient of Corollary \ref{cor:primitive-residue}.
Thus the primitive criterion detects uniform non-associativity that
Theorem \ref{theoremnonassociativeforall} does not detect. 
\end{rem}

\subsection{Controlling the correction by ideals}

The next statement keeps track of the coefficient ideals used in the
proof of Theorem~\ref{theoremnonassociativeforall}.
\begin{prop}[An ideal-valued criterion]
\label{prop:error-ideals} Let $a_{1},a_{2},a_{3}\in F$ be homogeneous,
and let $J,K,L\subseteq R$ be ideals satisfying 
\[
\mathrm{d}F\subseteq JF,\qquad\mathrm{d}a_{1},\mathrm{d}a_{2},\mathrm{d}a_{3}\in KF,\qquad a_{1}a_{2},a_{2}a_{3}\in LF.
\]
For every multiplication $\mu_{h}$ on $F$, 
\[
[a_{1},a_{2},a_{3}]_{\mu_{h}}-[a_{1},a_{2},a_{3}]_{\mu}\in(K+JL+J^{2})F+(\mathrm{d}F)a_{3}+a_{1}\mathrm{d}F.
\]
All products on the right are taken using $\mu$. Consequently, if
the reference associator lies outside this submodule, every such $\mu_{h}$
is non-associative at the same triple. 
\end{prop}

\begin{proof}
Let $\mu_{h}$ be a competing multiplication, and choose a homotopy
$h$ with $\mu_{h}=\mu+\mathrm{d}h+h\mathrm{d}$. Put $t=a_{1}\otimes a_{2}\otimes a_{3}$.
By \eqref{eq:homotopy-associator},
\[
[t]_{\mu_{h}}-[t]_{\mu}=(\mathrm{d}H+H\mathrm{d})(t),
\]
where 
\[
H=[\cdot]_{\mu,h}+[\cdot]_{h,\mu_{h}}=[\cdot]_{\mu,h}+[\cdot]_{h,\mu}+[\cdot]_{h,\mathrm{d}h}+[\cdot]_{h,h\mathrm{d}}.
\]

The hypotheses give 
\[
H\mathrm{d}(t)\in KF,\qquad\mathrm{d}[t]_{h,\mu}\in JLF,\qquad\mathrm{d}[t]_{h,\mathrm{d}h},\ \mathrm{d}[t]_{h,h\mathrm{d}}\in J^{2}F.
\]
The first inclusion uses $\mathrm{d}a_{i}\in KF$. For the second,
the adjacent products lie in $LF$, and applying $\mathrm{d}$ contributes
a factor of $J$. For the last two, both $\mathrm{d}h$ and $h\mathrm{d}$
have image in $JF$; applying the outer $h$ preserves this coefficient
ideal, and applying $\mathrm{d}$ contributes another factor of $J$.

It remains to consider $\mathrm{d}[t]_{\mu,h}$. Set $u=h(a_{1},a_{2})$
and $v=h(a_{2},a_{3})$. The Koszul sign rule and the Leibniz identity
give
\[
\begin{aligned}\mathrm{d}[t]_{\mu,h} & =\mathrm{d}\bigl(ua_{3}-(-1)^{|a_{1}|}a_{1}v\bigr)\\
 & \equiv(\mathrm{d}u)a_{3}-a_{1}\mathrm{d}v\pmod{KF}.
\end{aligned}
\]
The omitted terms contain $\mathrm{d}a_{1}$ or $\mathrm{d}a_{3}$,
so they belong to $KF$. The remaining terms lie in $(\mathrm{d}F)a_{3}+a_{1}\mathrm{d}F$.
Combining these inclusions proves the assertion.
\end{proof}
For a minimal resolution over a local ring, taking $J=L=\mathfrak{m}$
and $K=\mathfrak{m}^{2}$ recovers Theorem~\ref{theoremnonassociativeforall},
using \eqref{eq:graded-opposite-order} to move $a_{3}$ to the left
modulo $\mathfrak{m}^{2}F$. More generally, for integers $s\ge1$
and $t,u\ge0$, if 
\[
\mathrm{d}F\subseteq\mathfrak{m}^{s}F,\qquad\mathrm{d}a_{i}\in\mathfrak{m}^{t}F,\qquad a_{1}a_{2},a_{2}a_{3}\in\mathfrak{m}^{u}F,
\]
then the coefficient ideal in the proposition is $\mathfrak{m}^{q}$,
where 
\[
q=\min\{t,s+u,2s\}.
\]
Thus the exponent two in the original criterion comes from minimality
and the coefficient orders of the correction terms. The identity $\mathrm{d}^{2}=0$
does not make an expression such as $\mathrm{d}h\mathrm{d}$ vanish.

For degree-one triples, Theorem \ref{thm:primitive-criterion} also
subsumes Proposition \ref{prop:error-ideals}, using coefficient ideal
$J+L$. Indeed, for
\[
N^{\prime}=(J+L)F_{4}+\mathrm{d}F_{5}+F_{3}a_{3}+a_{1}F_{3},
\]
the Leibniz rule gives
\[
\mathrm{d}N^{\prime}\subseteq(K+JL+J^{2})F_{3}+(\mathrm{d}F_{3})a_{3}+a_{1}\mathrm{d}F_{3}.
\]
Thus the associator nonmembership condition forces its primitive outside
$N^{\prime}$. The two associator criteria retain their scope in arbitrary
homogeneous degrees.
\begin{example}[Powers of Avramov's example]
\label{ex:powered-avramov}Let $R=\Bbbk[x,y,z,w]$, let $\mathfrak{m}=\langle x,y,z,w\rangle$,
and put
\[
I_{N}=\langle x^{2N},w^{2N},z^{N}w^{N},x^{N}y^{N},y^{N}z^{N}\rangle,\qquad N\ge2.
\]
Apply the substitution $(x,y,z,w)\mapsto(x^{N},y^{N},z^{N},w^{N})$
to the resolution and reference multiplication in Subsection \ref{sec:avramov-uniform}.
The ring $R$ is free over $\Bbbk[x^{N},y^{N},z^{N},w^{N}]$, with
basis the monomials whose exponents are less than $N$. Thus substitution
preserves exactness. All differential entries remain in $\mathfrak{m}$,
so the resulting resolution $F$ is minimal. Retain the earlier basis
labels and write $F_{4}=Re_{12345}$.

The substituted associator is
\begin{equation}
[e_{1},e_{5},e_{2}]_{\mu}=y^{N}e_{123}-z^{N}e_{124}+w^{N}e_{1345}-x^{N}e_{2345}=-\mathrm{d}e_{12345}.\label{eq:powered-avramov-associator}
\end{equation}
It lies in $\mathfrak{m}^{N}F_{3}$, so Theorem \ref{theoremnonassociativeforall}
does not detect it. Put $\mathfrak{q}_{N}=\langle x^{N},y^{N},z^{N},w^{N}\rangle$.
The substituted formulas give
\[
\mathrm{d}F\subseteq\mathfrak{q}_{N}F,\qquad\mathrm{d}e_{1},\mathrm{d}e_{5},\mathrm{d}e_{2}\in\mathfrak{q}_{N}^{2},\qquad e_{1}e_{5},e_{5}e_{2}\in\mathfrak{q}_{N}F_{2}.
\]
The only nonzero products of $e_{1}$ or $e_{2}$ with the basis of
$F_{3}$ are
\[
e_{1}e_{2345}=x^{N}e_{12345},\qquad e_{2}e_{1345}=-w^{N}e_{12345}.
\]
The Leibniz rule therefore puts the endpoint terms $(\mathrm{d}F_{3})e_{2}+e_{1}\mathrm{d}F_{3}$
in $\mathfrak{q}_{N}^{2}F_{3}$. Applying Proposition \ref{prop:error-ideals}
with $J=L=\mathfrak{q}_{N}$ and $K=\mathfrak{q}_{N}^{2}$ gives
\[
[e_{1},e_{5},e_{2}]_{\mu_{h}}\equiv-\mathrm{d}e_{12345}\not\equiv0\pmod{\mathfrak{q}_{N}^{2}F_{3}}
\]
for every multiplication $\mu_{h}$.

Corollary \ref{cor:primitive-residue} gives the same conclusion after
localization at $\mathfrak{m}$. A primitive is $\tau=-e_{12345}$.
Its residue is nonzero, while the endpoint products vanish modulo
$\mathfrak{m}$. For $N=2$, this gives a quadratic associator detected
by the primitive and ideal-valued criteria. In general, the reference
associator has coefficient order exactly $N$. The conclusions hold
over every field, and after localization or completion at $\mathfrak{m}$.
\end{example}

\section{Permanence and examples}

\label{sec:examples}

In this section, we consider monomial ideals in polynomial rings over
a field $\Bbbk$ and their minimal free resolutions over the respective
polynomial rings, equipped with the fine multigrading and with degree-zero
module equal to the base ring.

\subsection{Passing the criteria from an elimination ideal}

The zero-new-variable part of the resolution below is the pruned resolution
of Boocher \cite[Theorem~2.3]{Boo12}. We include the block-matrix
argument to record the specialization retraction and the compatibility
with the reference multiplication needed for our obstruction. Geller,
Martin, and Potts-Rubin prove that commutative DG algebra structures
descend under this pruning for squarefree monomial ideals \cite[Theorem~5.7 and Corollary~5.8]{GMP26}.
Our permanence statement instead carries a specified fixed-triple
obstruction from the smaller resolution to the larger one, and allows
arbitrary multiplications as competitors.

Let $R=\Bbbk[x_{1},\ldots,x_{n}]$, let $S=R[y_{1},\ldots,y_{m}]$,
and write 
\[
\mathfrak{m}=\langle x_{1},\ldots,x_{n}\rangle,\qquad\mathfrak{n}=\langle x_{1},\ldots,x_{n},y_{1},\ldots,y_{m}\rangle.
\]
Let $J\subsetneq S$ be a monomial ideal, and let $I=J\cap R$ be
its \emph{elimination ideal}, obtained by eliminating $y_{1},\ldots,y_{m}$.
Equivalently, $I$ is the contraction of $J$ along $R\hookrightarrow S$.
In this monomial setting, elimination also agrees with specialization
to zero: if $\pi\colon S\to R$ sets $y_{1}=\cdots=y_{m}=0$, then
$\pi(J)=I$. Indeed, each monomial generator either involves a new
variable and vanishes, or belongs to $R$ and survives unchanged.
For squarefree ideals, the Stanley--Reisner complex of $I$ is the
induced subcomplex of the Stanley--Reisner complex of $J$ on the
original vertices.
\begin{lem}
\label{lem:specialization-retraction} Let $J\subsetneq S$ be a monomial
ideal, put $I=J\cap R$, and let $F$ and $G$ be the minimal free
resolutions of $R/I$ over $R$ and of $S/J$ over $S$, respectively.
After a multigraded identification, $F\otimes_{R}S$ is a subcomplex
of $G$, and 
\[
G\otimes_{S}R\cong F\oplus Q
\]
as $R$-complexes, where $S\to R$ sets $y_{1},\ldots,y_{m}$ equal
to zero. Every multigraded multiplication $\mu$ on $F$ extends to
a multigraded multiplication $\mu'$ on $G$. If $\mu$ is commutative,
the extension $\mu'$ can be chosen commutative. For this extension,
there is an $R$-linear map $\rho\colon G\to F$ satisfying 
\[
\rho(sa)=s|_{y=0}\,\rho(a),\qquad\rho|_{F}=\operatorname{id}_{F},\qquad\rho\mathrm{d}=\mathrm{d}\rho,\qquad\rho\mu'(a,b)=\mu(\rho(a),\rho(b)).
\]
\end{lem}

\begin{proof}
Choose a multihomogeneous basis of $G$. Its basis degrees have nonnegative
coordinates. Let $E_{i}$ be the $S$-span of the basis elements of
$G_{i}$ whose degrees in all the new variables are zero, and let
$P_{i}$ be the $S$-span of the remaining basis elements. Thus $G_{i}=E_{i}\oplus P_{i}$
as free $S$-modules. A differential term from $E_{i}$ to $P_{i-1}$
would require a coefficient with a negative degree in some new variable,
which is impossible. Hence $E$ is a subcomplex, and its differential
has coefficients in $R$.

Taking the part of new-variable degree zero in the augmented resolution
$G\to S/J$ over $S$ preserves exactness. This part is the $R$-span
of the basis elements defining $E$, and its zeroth homology is $R/I$,
since $J\cap R=I$. It is a minimal free resolution over $R$, which
we identify with $F$. Consequently $E=F\otimes_{R}S$.

To see the splitting after specialization explicitly, write the differential
with respect to $G_{i}=E_{i}\oplus P_{i}$ as 
\[
\mathrm{d}_{i}^{G}=\begin{pmatrix}\mathrm{d}_{i}^{E} & A_{i}\\
0 & B_{i}
\end{pmatrix}.
\]
A term from a basis element of $P_{i}$ to one of $E_{i-1}$ has a
coefficient of positive degree in at least one new variable. Thus
every entry of $A_{i}$ belongs to $\langle y_{1},\ldots,y_{m}\rangle$
and vanishes after setting $y=0$. Since $E\otimes_{S}R\cong F$,
the specialized differential is 
\[
\mathrm{d}_{i}^{G\otimes_{S}R}=\begin{pmatrix}\mathrm{d}_{i}^{F} & 0\\
0 & \overline{B}_{i}
\end{pmatrix},\qquad\overline{B}_{i}=B_{i}|_{y=0}.
\]
In particular, the modules $Q_{i}=P_{i}\otimes_{S}R$ with differentials
$\overline{B}_{i}$ form an $R$-complex $Q$, and the displayed block
diagonal differential gives 
\[
G\otimes_{S}R\cong F\oplus Q
\]
as $R$-complexes. The modules $P_{i}$ need not form a subcomplex
of $G$ before specialization. Define $\rho$ as the composite 
\[
G\longrightarrow G\otimes_{S}R\cong F\oplus Q\longrightarrow F,
\]
where the first map specializes coefficients and the last map is the
projection. This proves the first three identities.

To extend an arbitrary multiplication, take $C=G\otimes_{S}G$. If
$\mu$ is commutative and we want a commutative extension, instead
take the quotient complex 
\[
C=(G\otimes_{S}G)/\langle u\otimes v-(-1)^{|u||v|}v\otimes u,\ u\otimes u\text{ for }|u|\text{ odd}\rangle,
\]
where the relations are taken on homogeneous elements. The differential
descends to this quotient. In particular, the differential of an odd
square is $\mathrm{d}u\otimes u-u\otimes\mathrm{d}u$, which is zero
modulo the commutativity relations. The complex $C$ is degreewise
free on ordered basis pairs in the first case, and on unordered basis
pairs, omitting repeated odd basis elements, in the second. In either
case, the unit pairs and the pairs from $E$ span a subcomplex which
is a direct summand in each degree. Prescribe the product there using
the unit and $\mu\otimes_{R}S$. These prescriptions agree on their
intersection.

Extend this chain map inductively to the remaining basis pairs. Each
such pair has total homological degree at least two. At each stage
its required differential is a cycle in a positive homological degree
of $G$, so exactness supplies a lift. We may choose the lift in the
required multidegree. The resulting map $C\to G$ gives a multiplication
$\mu'$ extending $\mu$, commutative in the second case. This construction
involves no division by two.

Finally, nonnegativity of the new-variable degrees shows that a product
with an input of positive new-variable degree also has positive new-variable
degree. The map $\rho$ kills all such elements. On inputs of zero
new-variable degree the product is $\mu$. This gives the last identity. 
\end{proof}
The splitting in Lemma~\ref{lem:specialization-retraction} is a
splitting after specialization. In general $F\otimes_{R}S$ is not
a direct summand of $G$ as an $S$-complex. For example, for $I=\langle x\rangle\subseteq\Bbbk[x]$
and $J=\langle x,y\rangle\subseteq\Bbbk[x,y]$, a chain retraction
from the Koszul resolution over $\Bbbk[x,y]$ would have to send the
generator with differential $y$ to $ae_{x}$ with $ax=y$, which
is impossible.
\begin{prop}[Persistence from an elimination ideal]
\label{prop:elimination-persistence} With the notation of Lemma~\ref{lem:specialization-retraction},
suppose a multigraded multiplication $\mu$ on $F$ and homologically
homogeneous elements $a_{1},a_{2},a_{3}\in F$ satisfy 
\[
\mathrm{d}a_{i}\in\mathfrak{m}^{2}F\quad(1\le i\le3),\qquad a_{1}a_{2},a_{2}a_{3}\in\mathfrak{m}F,
\]
and 
\[
[a_{1},a_{2},a_{3}]_{\mu}\notin\mathfrak{m}^{2}F+a_{1}\mathrm{d}F+a_{3}\mathrm{d}F.
\]
Then every multiplication on $G$ is non-associative at the same triple
$(a_{1},a_{2},a_{3})$, viewed in $G$. The competing multiplication
need not respect the internal grading or be graded-commutative. 
\end{prop}

\begin{proof}
Choose $\mu'$ and $\rho$ as in Lemma~\ref{lem:specialization-retraction}.
The differential and the reference multiplication on the embedded
$F$ are unchanged, so the first two hypotheses persist with $\mathfrak{n}$
in place of $\mathfrak{m}$. Moreover, 
\[
\rho(\mathfrak{n}^{2}G)\subseteq\mathfrak{m}^{2}F,\qquad\rho(a_{i}\mathrm{d}G)\subseteq a_{i}\mathrm{d}F\quad(i=1,3).
\]
If the associator belonged to $\mathfrak{n}^{2}G+a_{1}\mathrm{d}G+a_{3}\mathrm{d}G$,
applying $\rho$ would contradict the assumed nonmembership in $F$.
Thus all three hypotheses of the uniform criterion hold for $\mu'$.

The quotient by $\mathfrak{n}^{2}G+a_{1}\mathrm{d}G+a_{3}\mathrm{d}G$
is annihilated by $\mathfrak{n}^{2}$, so this nonzero class remains
nonzero after localization at $\mathfrak{n}$. Apply Theorem~\ref{theoremnonassociativeforall}
and Remark~\ref{rem:scope-uniform}. Every polynomial multiplication
localizes, which proves the assertion. 
\end{proof}
The same specialization map preserves the residue-field primitive
obstruction of Theorem~\ref{thm:primitive-criterion}; we record
this stronger permanence statement in Corollary~\ref{cor:primitive-elimination}.
\begin{cor}[Persistence of the primitive obstruction]
\label{cor:primitive-elimination} Use the notation of Lemma~\ref{lem:specialization-retraction},
and write $\mathfrak{m}=\langle x_{1},\ldots,x_{n}\rangle$ and $\mathfrak{n}=\langle x_{1},\ldots,x_{n},y_{1},\ldots,y_{m}\rangle$.
Let $\mu$ be a multigraded multiplication on $F$, let $a_{1},a_{2},a_{3}\in F_{1}$,
and suppose 
\[
a_{1}a_{2},a_{2}a_{3}\in\mathfrak{m}F_{2},\qquad\tau_{\mu}\notin\mathfrak{m}F_{4}+\mathrm{d}F_{5}+a_{3}F_{3}+a_{1}F_{3},
\]
where $\mathrm{d}\tau_{\mu}=[a_{1},a_{2},a_{3}]_{\mu}$. After identifying
$F$ with its old-variable part in $G$, every multiplication on $G$
is non-associative at $(a_{1},a_{2},a_{3})$. 
\end{cor}

\begin{proof}
Choose the extended reference multiplication $\mu'$ and the specialization
map $\rho$ from Lemma~\ref{lem:specialization-retraction}. The
element $\tau_{\mu}$ is still a primitive of the reference associator
in $G$, because both the differential and the reference product restrict
to those on $F$. By \eqref{eq:graded-opposite-order}, the submodules
$F_{3}a_{3}+\mathfrak{m}F_{4}$ and $a_{3}F_{3}+\mathfrak{m}F_{4}$
coincide; the same holds for $G$ and $\mathfrak{n}$. Minimality
gives $\mathrm{d}G_{3}\subseteq\mathfrak{n}G_{2}$, and the two adjacent
products belong to $\mathfrak{n}G_{2}$. Moreover, 
\[
\begin{split}\rho\bigl(\mathfrak{n}G_{4}+\mathrm{d}G_{5}+a_{3}G_{3}+a_{1}G_{3}\bigr)\subseteq{} & \mathfrak{m}F_{4}+\mathrm{d}F_{5}\\
 & +a_{3}F_{3}+a_{1}F_{3}.
\end{split}
\]
Here we use semilinearity, $\rho\mathrm{d}=\mathrm{d}\rho$, multiplicativity
of $\rho$ for the reference products, and $\rho(a_{i})=a_{i}$. Since
$\rho(\tau_{\mu})=\tau_{\mu}$, the assumed nonmembership in $F$
implies the corresponding nonmembership in $G$. Theorem~\ref{thm:primitive-criterion},
with coefficient ideal $\mathfrak{n}$, proves the result. No restriction
of an arbitrary competing product on $G$ to $F$ is used. 
\end{proof}

\subsection{The path on six vertices}

\label{subsec:path}

Let 
\[
R=\Bbbk[x_{1},\ldots,x_{6}],\qquad I=\langle x_{1}x_{2},x_{2}x_{3},x_{3}x_{4},x_{4}x_{5},x_{5}x_{6}\rangle.
\]
This is the edge ideal of the path on six vertices, and it is the
polarization of Avramov's ideal. Its nonexistence of a graded-commutative
DG algebra structure is also discussed in \cite[Example~3.8]{Kat19}.
We record the calculation needed to apply our uniform criterion.

Use the same cellular basis as in Subsection~\ref{sec:avramov-uniform},
now denoted by $\eta_{\sigma}$, with vertex labels ordered as 
\[
\mathrm{d}\eta_{1}=x_{1}x_{2},\quad\mathrm{d}\eta_{2}=x_{5}x_{6},\quad\mathrm{d}\eta_{3}=x_{4}x_{5},\quad\mathrm{d}\eta_{4}=x_{2}x_{3},\quad\mathrm{d}\eta_{5}=x_{3}x_{4}.
\]
The labels are obtained by polarizing $x^{2}$ to $x_{1}x_{2}$ and
$w^{2}$ to $x_{5}x_{6}$, with $x,y,z,w$ replaced by $x_{2},x_{3},x_{4},x_{5}$,
respectively. Polarization preserves the divisibility subcomplexes
in the cellular acyclicity criterion, so the resulting labeled complex
is a minimal resolution over $R$, with ranks $(1,5,7,4,1)$. Set
$\tau_{\mathrm{P}}=\eta_{12345}$. The top differential is 
\begin{equation}
\mathrm{d}\tau_{\mathrm{P}}=-x_{3}\eta_{123}+x_{4}\eta_{124}-x_{6}\eta_{1345}+x_{1}\eta_{2345}.\label{eq:path-top}
\end{equation}
Choose a multigraded commutative multiplication $\mu_{\mathrm{P}}$.
As before, multigrading and the Leibniz rule force 
\begin{align*}
\eta_{1}\eta_{5} & =x_{4}\eta_{14}+x_{1}\eta_{45}, & \eta_{2}\eta_{5} & =x_{3}\eta_{23}+x_{6}\eta_{35},\\
\eta_{2}\eta_{14} & =-\eta_{124}, & \eta_{1}\eta_{23} & =\eta_{123},\\
\eta_{2}\eta_{45} & =-\eta_{2345}, & \eta_{1}\eta_{35} & =\eta_{1345}.
\end{align*}
Hence 
\[
[\eta_{1},\eta_{5},\eta_{2}]_{\mu_{\mathrm{P}}}=-\mathrm{d}\tau_{\mathrm{P}}.
\]
The endpoint products with $F_{3}$ vanish except for 
\[
\eta_{1}\eta_{2345}=x_{2}\tau_{\mathrm{P}},\qquad\eta_{2}\eta_{1345}=-x_{5}\tau_{\mathrm{P}}.
\]
Since $\mathrm{d}\eta_{1},\mathrm{d}\eta_{2}\in\mathfrak{m}^{2}$,
the Leibniz rule gives 
\[
\eta_{1}\mathrm{d}F_{3}+\eta_{2}\mathrm{d}F_{3}\subseteq\mathfrak{m}^{2}F_{3}.
\]
Thus the nonzero linear associator in \eqref{eq:path-top} satisfies
the nonmembership condition. The other two hypotheses follow from
the quadratic vertex differentials and the displayed products. Theorem~\ref{theoremnonassociativeforall}
and Remark~\ref{rem:scope-uniform} therefore give uniform non-associativity
at $(\eta_{1},\eta_{5},\eta_{2})$ for every multiplication.

\subsection{A Gorenstein example}

\label{sec:gorenstein}

Denham and Suciu identify $\operatorname{Bier}(C_{4})$ with their
eight-vertex sphere and exhibit a nonzero triple Massey product in
the characteristic-zero setting of \cite[Section 8, especially Examples 8.1.1 and 8.6.1]{DS07}.
Their Remark 8.1.2 places the associated moment-angle manifold in
Baskakov's family \cite{Bas03}. Grade-five Gorenstein counterexamples
were previously constructed in \cite{Sri96}. The theorem below identifies
an explicit fixed triple at which every multiplication fails over
every field.

More explicitly, take $C_{4}$ on vertices $1,2,3,4$ with edges $12,23,34,41$,
and number the vertices of its Alexander dual by $i'=i+4$. The vertex
map 
\[
(1,2,3,4,5,6,7,8)\longmapsto(1,3,7,6,2,4,8,5)
\]
identifies the Stanley--Reisner complex of \eqref{eq:gorenstein-ideal}
with the deleted join $\operatorname{Bier}(C_{4})=C_{4}*_{\Delta}C_{4}^{\vee}$.
\begin{thm}
\label{thm:gorenstein} Let $S=\Bbbk[x_{1},\ldots,x_{8}]$ and let
\begin{equation}
\begin{split}J=\langle\, & x_{1}x_{2},x_{2}x_{3},x_{3}x_{4},x_{4}x_{5},x_{5}x_{6},\\
 & x_{3}x_{7},x_{6}x_{7},x_{1}x_{8},x_{4}x_{8},x_{7}x_{8}\rangle.
\end{split}
\label{eq:gorenstein-ideal}
\end{equation}
Then $S/J$ is Gorenstein of dimension $3$ and codimension $5$ over
every field. Let $G$ be its minimal multigraded free resolution over
$S$, and let $b_{i}\in G_{1}$ be the basis element with $\mathrm{d}b_{i}=x_{i}x_{i+1}$
for $1\le i\le5$. Every multiplication on $G$ is non-associative
at $(b_{1},b_{3},b_{5})$. In particular, $G$ admits no DG algebra
structure, even when compatibility with the internal grading and graded
commutativity are not required. 
\end{thm}

\begin{proof}
We first verify the Gorenstein assertion. Write $ijk$ for the face
$\{i,j,k\}$. The Stanley--Reisner complex of the path ideal $I$
has facets 
\[
135,\quad136,\quad146,\quad246,\quad25.
\]
Adjoining $x_{7}$ and the generators $x_{3}x_{7},x_{6}x_{7}$ gives
the complex $B$ with facets 
\[
135,\quad136,\quad146,\quad246,\quad147,\quad157,\quad247,\quad257.
\]
This is a triangulated disk. Indeed, in the indicated order each triangle
after the first is attached along one boundary edge or two consecutive
boundary edges; the successive attaching intersections are 
\[
13,\quad16,\quad46,\quad14,\quad15\cup17,\quad24\cup47,\quad27\cup57.
\]
The boundary of $B$ is the cycle $2-5-3-6-2$. Adjoining $x_{8}$
and the generators $x_{1}x_{8},x_{4}x_{8},x_{7}x_{8}$ adds exactly
the four triangles 
\[
258,\quad358,\quad368,\quad268.
\]
Thus the Stanley--Reisner complex of $J$ is $B\cup(8*\partial B)$,
a triangulated $2$-sphere. Its Stanley--Reisner ring is Gorenstein
over every field \cite[Theorem~II.5.1]{Sta96}. Since the facets have
three vertices and $S$ has eight variables, its dimension is $3$
and its codimension is $5$.

Finally, $J\cap\Bbbk[x_{1},\ldots,x_{6}]=I$. By Lemma~\ref{lem:specialization-retraction},
we may identify the path resolution with the part of $G$ of zero
degree in $x_{7},x_{8}$. Under this identification, $(\eta_{1},\eta_{5},\eta_{2})$
becomes $(b_{1},b_{3},b_{5})$. The calculation in Subsection~\ref{subsec:path}
and Proposition~\ref{prop:elimination-persistence} prove the asserted
uniform non-associativity. 
\end{proof}

\section{Katthän's example and relative non-associativity}

\label{sec:katthan}

Let $\Bbbk$ be a field, let $R=\Bbbk[x,y,z,w]$, let $\mathfrak{m}=\langle x,y,z,w\rangle$,
let $I=\langle x^{2},w^{2},zw,xy,y^{2}z^{2}\rangle$, and let $F$
be the minimal multigraded free resolution of $R/I$ over $R$, with
$F_{0}=R$. The ideal $I$ is obtained from Avramov's ideal by replacing
its last generator $yz$ by $y^{2}z^{2}$.

Katthän's theorem excludes graded-commutative DG algebra structures
whose multiplication respects the standard internal grading~\cite[Theorem~5.1]{Kat19}.
Taking the degree-zero component of a multiplication preserves the
chain-map and unit conditions, but need not preserve associativity,
since components of opposite internal degrees also contribute to the
degree-zero associativity equation. Removing the grading hypothesis
therefore requires an additional argument. Moreover, Example~\ref{ex:cancel}
shows that even among standard-graded commutative multiplications,
the associator at $(e_{1},e_{5},e_{2})$ can be cancelled. Our proof
instead controls several triples simultaneously and shows that no
multiplication on $F$ is associative, even without graded commutativity.

\subsection{The resolution}

In this subsection we describe $F$ as a resolution over $R$. The
ideal $I$ is strongly generic: no variable appears with the same
nonzero exponent in two distinct minimal generators. Its Scarf complex
is 
\[
\Delta=\langle\{1,2,3,4\},\{3,4,5\}\rangle.
\]
By \cite[Theorem~3.2]{BPS98}, this complex supports the minimal free
resolution over $R$; see also \cite[Chapter~6]{MS05}. One can visualize
$F$ as being supported on this simplicial complex below:
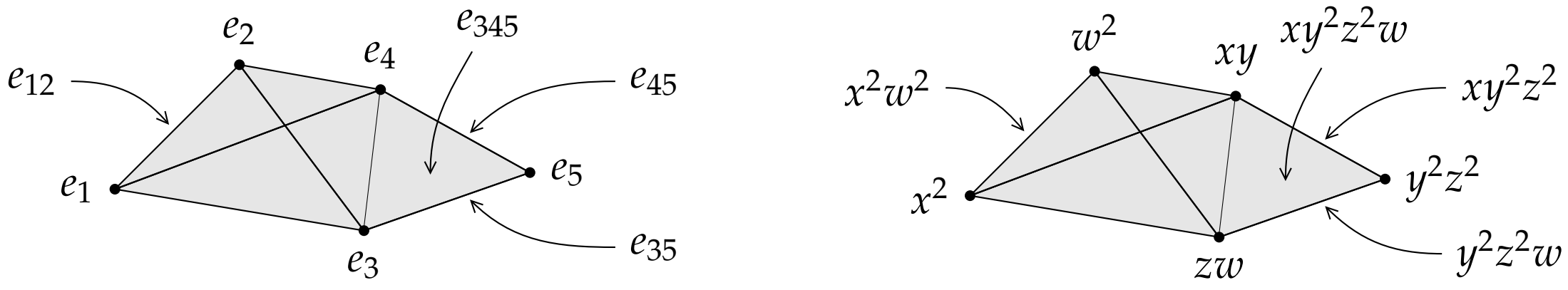
\begin{center}
\begin{tabular}{cccc}
\begin{tikzpicture}[scale=0.77]
\draw[fill=gray!20] (0,0) -- (3,-0.5) -- (3.2,1.2)-- (0,0); 
\draw[fill=gray!20] (0,0) -- (1.5,1.5) -- (3.2,1.2)-- (0,0); \draw[] (1.5,1.5) -- (3,-0.5); 
\draw[fill=gray!20] (3,-0.5) -- (5,0.2) -- (3.2,1.2); 
\draw[color=black!100] (1.5,1.5) -- (3,-0.5) -- (5,0.2);
\draw[color=black!100] (0,0) -- (3.2,1.2) -- (5,0.2);

\node[circle, fill=black, inner sep=1pt, label=left:$e_1 $] (a) at (0,0) {}; 
\node[circle, fill=black, inner sep=1pt, label=above:$e_2 $] (b) at (1.5,1.5) {}; 
\node[circle, fill=black, inner sep=1pt, label=below:$e_3 $] (c) at (3,-0.5) {};
\node[circle, fill=black, inner sep=1pt, label=above:$e_4 $] (d) at (3.2,1.2) {}; 
\node[circle, fill=black, inner sep=1pt, label=right:$e_5 $] (e) at (5,0.2) {};

\node[] (x) at (-1,1.3) {$ e_{12} $}; \node[] (x') at (0.8,0.6) {};
\draw[-{Straight Barb[length=3pt,width=3pt]}] (x) edge[out=0, in=130] node[below] {} (x');

\node[] (y) at (6.5,1.3) {$ e_{45} $}; \node[] (y') at (4.1,0.5) {}; 
\draw[-{Straight Barb[length=3pt,width=3pt]}] (y) edge[out=180, in=45] node[below] {} (y');

\node[] (z) at (6.5,-0.7) {$e_{35} $}; \node[] (z') at (4.1,0.05) {};
\draw[-{Straight Barb[length=3pt,width=3pt]}] (z) edge[out=180, in=-45] node[below] {} (z');

\node[] (w) at (4.5,2) {$e_{345} $}; \node[] (w') at (3.8,0) {};
\draw[-{Straight Barb[length=3pt,width=3pt]}] (w) edge[out=-120, in=90] node[below] {} (w');
\end{tikzpicture} &  &  & \begin{tikzpicture}[scale=0.77]

\draw[fill=gray!20] (0,0) -- (3,-0.5) -- (3.2,1.2)-- (0,0); 
\draw[fill=gray!20] (0,0) -- (1.5,1.5) -- (3.2,1.2)-- (0,0); 
\draw[] (1.5,1.5) -- (3,-0.5); \draw[fill=gray!20] (3,-0.5) -- (5,0.2) -- (3.2,1.2);
\draw[color=black!100] (1.5,1.5) -- (3,-0.5) -- (5,0.2);
\draw[color=black!100] (0,0) -- (3.2,1.2) -- (5,0.2);

\node[circle, fill=black, inner sep=1pt, label=left:$x^2 $] (a) at (0,0) {};
\node[circle, fill=black, inner sep=1pt, label=above:$w^2 $] (b) at (1.5,1.5) {};
\node[circle, fill=black, inner sep=1pt, label=below:$zw $] (c) at (3,-0.5) {};
\node[circle, fill=black, inner sep=1pt, label=above:$xy $] (d) at (3.2,1.2) {}; 
\node[circle, fill=black, inner sep=1pt, label=right:$y^2 z^2 $] (e) at (5,0.2) {};

\node[] (x) at (-1,1.3) {$ x^2 w^2 $}; \node[] (x') at (0.8,0.6) {};
\draw[-{Straight Barb[length=3pt,width=3pt]}] (x) edge[out=0, in=130] node[below] {} (x');

\node[] (y) at (6.5,1.3) {$ x y^2 z^2 $}; \node[] (y') at (4.1,0.5) {}; 
\draw[-{Straight Barb[length=3pt,width=3pt]}] (y) edge[out=180, in=45] node[below] {} (y');

\node[] (z) at (6.5,-0.7) {$ y^2 z^2 w $}; \node[] (z') at (4.1,0.05) {};
\draw[-{Straight Barb[length=3pt,width=3pt]}] (z) edge[out=180, in=-45] node[below] {} (z');

\node[] (w) at (4.5,2) {$ x y^2 z^2 w $}; \node[] (w') at (3.8,0) {}; 
\draw[-{Straight Barb[length=3pt,width=3pt]}] (w) edge[out=-120, in=90] node[below] {} (w');

\end{tikzpicture}\tabularnewline
\end{tabular}
\par\end{center}

Its free modules are 
\begin{align}
F_{0} & =R,\label{eq:modules}\nonumber \\
F_{1} & =Re_{1}+Re_{2}+Re_{3}+Re_{4}+Re_{5},\nonumber \\
F_{2} & =Re_{12}+Re_{13}+Re_{14}+Re_{23}+Re_{24}+Re_{34}+Re_{35}+Re_{45},\\
F_{3} & =Re_{123}+Re_{124}+Re_{134}+Re_{234}+Re_{345},\nonumber \\
F_{4} & =Re_{1234}.\nonumber 
\end{align}
For an ordered face $\sigma=(i_{1}<\cdots<i_{p})$, let $m_{\sigma}$
be the least common multiple of its vertex labels. Our orientation
convention is 
\begin{equation}
\mathrm{d}e_{\sigma}=\sum_{j=1}^{p}(-1)^{j-1}\frac{m_{\sigma}}{m_{\sigma\setminus\{i_{j}\}}}e_{\sigma\setminus\{i_{j}\}}.\label{eq:taylor-d}
\end{equation}
In particular, we have $\mathrm{d}e_{1}=x^{2}$, $\mathrm{d}e_{2}=w^{2}$,
$\mathrm{d}e_{3}=zw$, $\mathrm{d}e_{4}=xy$, and $\mathrm{d}e_{5}=y^{2}z^{2}$.
Similarly we have 
\begin{align}
\mathrm{d}e_{123} & =x^{2}e_{23}-we_{13}+ze_{12},\label{eq:d3}\nonumber \\
\mathrm{d}e_{124} & =xe_{24}-w^{2}e_{14}+ye_{12},\nonumber \\
\mathrm{d}e_{134} & =xe_{34}-zwe_{14}+ye_{13},\\
\mathrm{d}e_{234} & =we_{34}-ze_{24}+xye_{23},\nonumber \\
\mathrm{d}e_{345} & =we_{45}-xe_{35}+yze_{34},\nonumber 
\end{align}
and 
\begin{equation}
\mathrm{d}e_{1234}=-ye_{123}+ze_{124}-we_{134}+xe_{234}.\label{eq:dtau}
\end{equation}

For completeness, exactness can also be seen without invoking the
Scarf construction. If 
\[
J=\langle x^{2},w^{2},zw,xy\rangle,\qquad f=y^{2}z^{2},
\]
then $J:f=\langle x,w\rangle$. The Taylor complex on the first four
generators is a minimal resolution of $R/J$ over $R$, and the Koszul
complex over $R$ on $x,w$ resolves $R/\langle x,w\rangle$. The
mapping cone lifting 
\[
0\longrightarrow R/\langle x,w\rangle\xrightarrow{\ f\ }R/J\longrightarrow R/I\longrightarrow0
\]
is the complex \eqref{eq:modules}--\eqref{eq:dtau}, with the indicated
choices of signs. All differential entries belong to $\mathfrak{m}$,
so this resolution is minimal. A more general version of this construction
appears in Subsection~\ref{sec:family}.

\subsection{The multiplication}

\label{sec:reference}

In this subsection we construct a commutative multiplication on $F$
which will serve as a reference. Let $T$ be the Taylor resolution
over $R$ on all five generators, with its usual multiplication $\mu_{T}$
as recalled in \cite[Section~1.2]{Kat19}, and let $\iota\colon F\hookrightarrow T$
be the inclusion. We specify a chain retraction $\pi\colon T\to F$.
It is the identity on every basis face of $F$ and has the following
values on the remaining faces: 
\begin{equation}
\begin{aligned}\pi(e_{15}) & =yz^{2}e_{14}+xe_{45}, & \pi(e_{25}) & =y^{2}ze_{23}+we_{35},\\
\pi(e_{125}) & =yz^{2}e_{124}+xyze_{234}-xwe_{345}, & \pi(e_{135}) & =yze_{134}-xe_{345},\\
\pi(e_{145}) & =0, & \pi(e_{235}) & =0,\\
\pi(e_{245}) & =-yze_{234}+we_{345}, & \pi(e_{1235}) & =yze_{1234},\\
\pi(e_{1245}) & =0, & \pi(e_{1345}) & =0,\\
\pi(e_{2345}) & =0, & \pi(e_{12345}) & =0.
\end{aligned}
\label{eq:pi}
\end{equation}
The Taylor differential directly gives $\mathrm{d}\pi=\pi\mathrm{d}$.
Since $\pi\iota=\operatorname{id}_{F}$, the formula 
\begin{equation}
\mu=\pi\mu_{T}\iota^{\otimes2}\label{eq:reference}
\end{equation}
defines a commutative multiplication. This fixes the reference product
on all pairs, not only on the pairs needed below. We write $ab=\mu(a,b)$.
By Lemma~\ref{lem:normalized}, every competing multiplication has
the form $\mu_{h}=\mu+\mathrm{d}h+h\mathrm{d}$, with $h$ normalized.
For a fixed such $h$, write $a\star b=\mu_{h}(a,b)$, reserving juxtaposition
for the reference multiplication.

The subcomplex on the first four vertices is closed under $\mu$,
and its multiplication is exactly its Taylor multiplication. In particular,
$\mu$ is associative on all triples from $e_{1},e_{2},e_{3},e_{4}$.
Its degree-one products on those generators are 
\begin{equation}
e_{1}e_{2}=e_{12},\quad e_{1}e_{3}=e_{13},\quad e_{1}e_{4}=xe_{14},\quad e_{2}e_{3}=we_{23},\quad e_{2}e_{4}=e_{24},\quad e_{3}e_{4}=e_{34}.\label{eq:core-products}
\end{equation}
Other products used in the proof are 
\begin{equation}
\begin{aligned}e_{1}e_{5} & =yz^{2}e_{14}+xe_{45}, & e_{2}e_{5} & =y^{2}ze_{23}+we_{35},\\
e_{1}e_{23} & =e_{123}, & e_{2}e_{14} & =-e_{124},\\
e_{1}e_{35} & =yze_{134}-xe_{345}, & e_{2}e_{45} & =-yze_{234}+we_{345}.
\end{aligned}
\label{eq:necessary-products}
\end{equation}
They give 
\begin{align*}
[e_{1},e_{5},e_{2}]_{\mu} & =-yz^{2}e_{124}-xyze_{234}+xwe_{345}+y^{2}ze_{123}+wyze_{134}-wxe_{345}\\
 & =-yz\,\mathrm{d}e_{1234}.
\end{align*}
Finally, the products $e_{i}F_{3}\subseteq F_{4}$ for $i\le4$ are
recorded in the following table; each entry denotes the coefficient
of $e_{1234}$. 
\begin{equation}
\begin{array}{c|rrrrr}
 & e_{123} & e_{124} & e_{134} & e_{234} & e_{345}\\
\hline e_{1} & 0 & 0 & 0 & x & 0\\
e_{2} & 0 & 0 & -w & 0 & 0\\
e_{3} & 0 & w & 0 & 0 & 0\\
e_{4} & -x & 0 & 0 & 0 & 0
\end{array}\label{eq:top-table}
\end{equation}
Consequently, 
\begin{equation}
e_{i}F_{3}+F_{3}e_{i}\subseteq\langle x,w\rangle F_{4}\qquad(1\le i\le4).\label{eq:top-vanishing}
\end{equation}
This property belongs to the fixed reference product. No corresponding
property is assumed for an arbitrary multiplication $\mu_{h}$. The
computation above gives 
\begin{equation}
[e_{1},e_{5},e_{2}]_{\mu}=-yz\,\mathrm{d}e_{1234}\in\mathfrak{m}^{3}F_{3}.\label{eq:old-associator}
\end{equation}
Thus the nonmembership hypothesis of Theorem~\ref{theoremnonassociativeforall}
fails for this triple and this reference multiplication.
\begin{example}
\label{ex:cancel} Define the normalized graded $R$-linear map $h\colon F\otimes_{R}F\to F$
of homological degree one by 
\begin{equation}
h(e_{1},e_{5})=e_{345},\qquad h(e_{5},e_{1})=-e_{345},\qquad h(e_{2},e_{34})=e_{1234}=h(e_{34},e_{2}),\label{eq:explicit-h}
\end{equation}
and set it equal to zero on all other ordered pairs of basis elements.
Then $\mu_{h}=\mu+\mathrm{d}h+h\mathrm{d}$ is a commutative multiplication,
and the scalar variation formula in Lemma \ref{lem:variation} gives
\begin{equation}
[e_{1},e_{5},e_{2}]_{\mu_{h}}=0,\qquad[e_{2},e_{3},e_{4}]_{\mu_{h}}=-\mathrm{d}e_{1234}\ne0.\label{eq:cancellation}
\end{equation}
Thus $(e_{1},e_{5},e_{2})$ cannot give a uniform non-associativity
obstruction: the first associator can be cancelled, while another
one becomes nonzero.
\end{example}

The obstruction will instead be relative to associativity on certain
triples involving the first four degree-one basis elements. Their
unit-coefficient products in \eqref{eq:core-products} constrain the
same homotopy coefficients which could cancel \eqref{eq:old-associator}. 

\subsection{The relative obstruction}

\label{sec:proof} 

The proof uses associativity on triples from the first four generators
to control the corrections to the associator at $(e_{1},e_{5},e_{2})$.
We proceed in four steps.
\begin{enumerate}
\item Use $F_{4}=Re_{1234}$ and injectivity of $\mathrm{d}\colon F_{4}\to F_{3}$
to express each associator on three degree-one inputs as a coefficient
times $\mathrm{d}e_{1234}$.
\item Localize at $\mathfrak{m}$ and reduce modulo $\langle x,w\rangle$.
The reference products $e_{i}F_{3}+F_{3}e_{i}$ vanish for $i\le4$,
while the reductions of the four products $e_{1}e_{2},e_{1}e_{3},e_{2}e_{4},e_{3}e_{4}$
give a basis of the relevant submodule of the reduced $F_{2}$.
\item Impose associativity on the $32$ triples specified below. Their equations
and Nakayama's lemma force the relevant reduced homotopy maps to vanish
on that submodule.
\item Evaluate the remaining associator at $(e_{1},e_{5},e_{2})$. Its coefficient
is $-yz$ times a unit, so it cannot vanish.
\end{enumerate}
Throughout this subsection, let $U=\{1,4\}$, let $V=\{2,3\}$, and
define the set of $32$ ordered triples 
\begin{equation}
\mathcal{C}=(U\times U\times V)\cup(V\times V\times U)\cup(U\times V\times V)\cup(V\times U\times U).\label{eq:C}
\end{equation}
Equivalently, $\mathcal{C}$ consists of the triples whose first and
third indices lie in opposite sets $U,V$; the middle index is arbitrary
in $\{1,2,3,4\}$. Repeated indices are allowed. 
\begin{thm}
\label{thm:main} The minimal free resolution $F$ of $R/I$ over
$R$ admits no associative multiplication. The same assertion holds
for the corresponding minimal resolutions over $R_{\mathfrak{m}}$
and over $\Bbbk\llbracket x,y,z,w\rrbracket$. These assertions hold
in every characteristic. 
\end{thm}

We prove the following stronger statement.
\begin{thm}[Relative associativity obstruction]
\label{thm:relative} Let $\mu_{h}$ be any multiplication on $F$.
If 
\[
[e_{i},e_{j},e_{k}]_{\mu_{h}}=0\qquad\text{for every }(i,j,k)\in\mathcal{C},
\]
then 
\[
[e_{1},e_{5},e_{2}]_{\mu_{h}}\ne0.
\]
In particular, among the $33$ triples in $\mathcal{C}\cup\{(1,5,2)\}$,
at least one is non-associative for every multiplication. 
\end{thm}

Theorem~\ref{thm:main} follows immediately from Theorem~\ref{thm:relative}.
The relative theorem asserts that associativity on $\mathcal{C}$
prevents cancellation of \eqref{eq:old-associator}.

\subsubsection{Two lemmas on homotopy coefficients}

The reduction to scalar associators uses the rank-one module $F_{4}$.
 All statements in this section apply after localization or completion
as well.

By Lemma~\ref{lem:normalized}, write $\mu_{h}=\mu+\mathrm{d}h+h\mathrm{d}$
with $h$ normalized. Since $F_{\ge5}=0$, its only potentially nonzero
components on pairs of positive homological degrees are 
\[
h\colon F_{1}\otimes_{R}F_{1}\to F_{3},\qquad h\colon F_{1}\otimes_{R}F_{2}\to F_{4},\qquad h\colon F_{2}\otimes_{R}F_{1}\to F_{4}.
\]
No symmetry of these components is imposed.

Fix such a normalized $h$ and write 
\begin{equation}
h_{ij}=h(e_{i},e_{j})\in F_{3},\qquad h(e_{i},a)=\ell_{i}(a)e_{1234},\qquad h(a,e_{i})=\ell'_{i}(a)e_{1234}\quad(a\in F_{2}),\label{eq:h-parameters}
\end{equation}
where $\ell_{i},\ell'_{i}\colon F_{2}\to R$ are $R$-linear maps.
Since $\mathrm{d}e_{i},\mathrm{d}e_{j}\in R$ and $h$ vanishes when
either input is a scalar, the $h\mathrm{d}$ term vanishes on $e_{i}\otimes e_{j}$.
Therefore 
\begin{equation}
e_{i}\star e_{j}=e_{i}e_{j}+\mathrm{d}h_{ij}.\label{eq:P}
\end{equation}
In particular, neither $h_{ii}=0$ nor $h_{ji}=-h_{ij}$ is being
assumed.

For any multiplication $\nu$, the element $[e_{i},e_{j},e_{k}]_{\nu}$
is a cycle in $F_{3}$. Indeed, the Leibniz formula for the associator
reduces its differential to associators with a scalar input, and those
vanish by unitality and $R$-bilinearity. Exactness and injectivity
of $\mathrm{d}\colon F_{4}\to F_{3}$ therefore give unique $a_{ijk}^{\nu}\in R$
such that 
\begin{equation}
[e_{i},e_{j},e_{k}]_{\nu}=a_{ijk}^{\nu}\,\mathrm{d}e_{1234}.\label{eq:scalar-def}
\end{equation}

\begin{lem}[Scalar variation formula]
\label{lem:variation} For every ordered triple $i,j,k\in\{1,\ldots,5\}$,
one has 
\[
\bigl(a_{ijk}^{\mu_{h}}-a_{ijk}^{\mu}\bigr)e_{1234}=[e_{i},e_{j},e_{k}]_{\mu,h}+[e_{i},e_{j},e_{k}]_{h,\mu_{h}}.
\]
Equivalently,
\begin{equation}
\begin{aligned}a_{ijk}^{\mu_{h}}e_{1234} & =a_{ijk}^{\mu}e_{1234}+h_{ij}e_{k}+e_{i}h_{jk}\\
 & \quad+\bigl(\ell'_{k}(e_{i}\star e_{j})-\ell_{i}(e_{j}\star e_{k})\bigr)e_{1234}.
\end{aligned}
\label{eq:variation}
\end{equation}
This identity does not require graded commutativity of $\mu_{h}$. 
\end{lem}

\begin{proof}
Let $H=[\cdot]_{\mu,h}+[\cdot]_{h,\mu_{h}}$ as in \eqref{eq:homotopy-associator}.
Normalization of $h$ and unitality of the multiplications imply that
the two generalized associators cancel whenever an input is a scalar.
Thus $H$ vanishes on such tensors. Since $\mathrm{d}e_{i},\mathrm{d}e_{j},\mathrm{d}e_{k}\in R$,
we have $H\mathrm{d}(e_{i}\otimes e_{j}\otimes e_{k})=0$. Consequently,
\eqref{eq:homotopy-associator} identifies the difference of the two
associators with $\mathrm{d}H(e_{i},e_{j},e_{k})$. Together with
\eqref{eq:scalar-def} and the injectivity of $\mathrm{d}\colon F_{4}\to F_{3}$,
this proves the first identity. Since $|e_{i}|=1$, evaluating the
generalized associators gives
\[
\begin{aligned}{}[e_{i},e_{j},e_{k}]_{\mu,h} & =h_{ij}e_{k}+e_{i}h_{jk},\\{}
[e_{i},e_{j},e_{k}]_{h,\mu_{h}} & =\bigl(\ell'_{k}(e_{i}\star e_{j})-\ell_{i}(e_{j}\star e_{k})\bigr)e_{1234}.
\end{aligned}
\]
Substituting these expressions yields \eqref{eq:variation}, the scalar
expression of \eqref{eq:homotopy-associator}.
\end{proof}
\Needspace{10\baselineskip} 

The following lemma controls the relevant homotopy coefficients by
Nakayama's lemma.
\begin{lem}
\label{lem:nakayama} Let $(B,\mathfrak{n})$ be a commutative local
ring, and let $G$ be a free $B$-module with basis $e_{12},e_{13},e_{24},e_{34}$.
Use $U,V,\mathcal{C}$ from \eqref{eq:C}. Let $p_{ij}\in G$ for
$i,j\in\{1,2,3,4\}$ be alternating, meaning $p_{ii}=0$ and $p_{ji}=-p_{ij}$
for all $i,j$, with 
\[
p_{12}=e_{12},\quad p_{13}=e_{13},\quad p_{24}=e_{24},\quad p_{34}=e_{34},\qquad p_{14}=p_{23}=0.
\]
Suppose $q_{ij}\in G$ satisfy $q_{ij}-p_{ij}\in\mathfrak{n}G$ for
all $i,j\in\{1,2,3,4\}$. Note that no alternation is assumed for
the $q_{ij}$. Let $\ell_{i},\ell'_{i}\in\operatorname{Hom}_{B}(G,B)$
for $i\in\{1,2,3,4\}$. If 
\begin{equation}
\ell'_{k}(q_{ij})=\ell_{i}(q_{jk})\qquad((i,j,k)\in\mathcal{C}),\label{eq:nak-equations}
\end{equation}
then all eight linear maps $\ell_{i},\ell'_{i}$ are zero. 
\end{lem}

\begin{proof}
Let $J\subseteq B$ be the ideal generated by the $32$ values of
the eight linear maps on the four displayed basis vectors, thus 
\[
J=\left\langle \begin{array}{l}
\ell_{i}(e_{12}),\ell_{i}(e_{13}),\ell_{i}(e_{24}),\ell_{i}(e_{34}),\\
\ell'_{i}(e_{12}),\ell'_{i}(e_{13}),\ell'_{i}(e_{24}),\ell'_{i}(e_{34})
\end{array}\;\middle|\;1\le i\le4\right\rangle _{B}.
\]

We show $J\subseteq\mathfrak{n}J$.

Fix $i\in U$. Let $j$ range through $U$ and $k$ through $V$.
Then $p_{ij}=0$, so $q_{ij}\in\mathfrak{n}G$, while the four elements
$p_{jk}$ run through a basis of $G$ up to signs. From \eqref{eq:nak-equations},
\[
\ell_{i}(p_{jk})=\ell'_{k}(q_{ij})-\ell_{i}(q_{jk}-p_{jk})\in\mathfrak{n}J.
\]
Thus $\ell_{i}(G)\subseteq\mathfrak{n}J$. The same argument with
$U,V$ interchanged applies to $i\in V$.

Next fix $k\in U$, let $j$ range through $U$, and let $i$ range
through $V$. Now $q_{jk}\in\mathfrak{n}G$, and the four $p_{ij}$
again form a basis up to signs. Thus 
\[
\ell'_{k}(p_{ij})=\ell_{i}(q_{jk})-\ell'_{k}(q_{ij}-p_{ij})\in\mathfrak{n}J.
\]
Interchanging $U,V$ gives the assertion for $k\in V$. Every generator
of $J$ belongs to $\mathfrak{n}J$. Hence $J=\mathfrak{n}J$, and
Nakayama's lemma gives $J=0$. 
\end{proof}

\subsubsection{Proof of the relative obstruction}
\begin{proof}[Proof of Theorem~\ref{thm:relative}]
It suffices to work over $R_{\mathfrak{m}}$: any multiplication
over $R$ localizes to one over $R_{\mathfrak{m}}$. Put 
\[
B=R_{\mathfrak{m}}/\langle x,w\rangle\cong\Bbbk[y,z]_{\langle y,z\rangle},\qquad\mathfrak{n}=\langle y,z\rangle B,
\]
and use a bar for reduction modulo $\langle x,w\rangle$. Inside $\overline{F}_{2}=F_{2}\otimes B$,
let 
\begin{equation}
G=B\bar{e}_{12}+B\bar{e}_{13}+B\bar{e}_{24}+B\bar{e}_{34}.\label{eq:W}
\end{equation}
Equations \eqref{eq:d3} give 
\begin{equation}
\begin{array}{c|ccccc}
\bar{a} & \bar{e}_{123} & \bar{e}_{124} & \bar{e}_{134} & \bar{e}_{234} & \bar{e}_{345}\\
\hline \mathrm{d}\bar{a} & z\bar{e}_{12} & y\bar{e}_{12} & y\bar{e}_{13} & -z\bar{e}_{24} & yz\bar{e}_{34}.
\end{array}\label{eq:reduced-d}
\end{equation}
Consequently $\mathrm{d}\overline{F}_{3}\subseteq\mathfrak{n}G$.
Put $p_{ij}=\overline{e_{i}e_{j}}$ and $q_{ij}=\overline{e_{i}\star e_{j}}$
for $i,j\le4$. Thus $p_{ij}$ is the reduction of the reference product,
whereas $q_{ij}$ is the reduced perturbed product. For $i,j\le4$,
\eqref{eq:core-products} and \eqref{eq:P} show 
\begin{equation}
q_{ij}\in G,\qquad q_{ij}-p_{ij}\in\mathfrak{n}G.\label{eq:reduced-P}
\end{equation}
The reference products satisfy $p_{14}=p_{23}=0$, and $p_{12},p_{13},p_{24},p_{34}$
form the displayed basis of $G$. These are the elements in Lemma~\ref{lem:nakayama},
with its abstract basis elements $e_{12},e_{13},e_{24},e_{34}$ identified
with $\bar{e}_{12},\bar{e}_{13},\bar{e}_{24},\bar{e}_{34}$, respectively.

The hypothesis of the theorem says $a_{ijk}^{\mu_{h}}=0$ for $(i,j,k)\in\mathcal{C}$.
The reference product is associative on the first four generators,
so $a_{ijk}^{\mu}=0$ there as well. The two $\mu$-terms in \eqref{eq:variation}
vanish after reduction by \eqref{eq:top-vanishing}. We therefore
obtain 
\[
\bar{\ell}'_{k}(q_{ij})=\bar{\ell}_{i}(q_{jk})\qquad((i,j,k)\in\mathcal{C}).
\]
Apply Lemma~\ref{lem:nakayama} with $q_{ij}=\overline{e_{i}\star e_{j}}$
and the $B$-linear maps $\bar{\ell}_{i}|_{G},\bar{\ell}'_{i}|_{G}$.
It gives 
\begin{equation}
\bar{\ell}_{i}(G)=0=\bar{\ell}'_{i}(G)\qquad(1\le i\le4).\label{eq:vanish-W}
\end{equation}

Now apply \eqref{eq:variation} to $(i,j,k)=(1,5,2)$. By \eqref{eq:top-vanishing},
the reference-product terms $h_{15}e_{2}+e_{1}h_{52}$ vanish modulo
$\langle x,w\rangle$. Equations \eqref{eq:P} and \eqref{eq:necessary-products}
give
\[
\overline{e_{1}\star e_{5}}=yz^{2}\bar{e}_{14}+\mathrm{d}\bar{h}_{15},\qquad\overline{e_{5}\star e_{2}}=-y^{2}z\bar{e}_{23}+\mathrm{d}\bar{h}_{52}.
\]
Since $h_{15},h_{52}\in F_{3}$, the differential table \eqref{eq:reduced-d}
places both $\mathrm{d}\bar{h}_{15}$ and $\mathrm{d}\bar{h}_{52}$
in $\mathrm{d}\bar{F}_{3}\subseteq\mathfrak{n}G$. Thus \eqref{eq:vanish-W}
gives $\bar{\ell}'_{2}(\mathrm{d}\bar{h}_{15})=0$ and $\bar{\ell}_{1}(\mathrm{d}\bar{h}_{52})=0$.
Using $a_{152}^{\mu}=-yz$, the scalar variation formula becomes
\begin{align}
\bar{a}_{152}^{\mu_{h}} & =-yz+\bar{\ell}'_{2}(\overline{e_{1}\star e_{5}})-\bar{\ell}_{1}(\overline{e_{5}\star e_{2}})\nonumber \\
 & =-yz+yz^{2}\bar{\ell}'_{2}(\bar{e}_{14})+y^{2}z\bar{\ell}_{1}(\bar{e}_{23})\nonumber \\
 & =-yz\bigl(1-z\bar{\ell}'_{2}(\bar{e}_{14})-y\bar{\ell}_{1}(\bar{e}_{23})\bigr).\label{eq:unit-obstruction}
\end{align}
The parenthesized factor is a unit of $B$, since its residue in $B/\mathfrak{n}$
is $1$. The ring $B$ is a domain and $yz\ne0$, so \eqref{eq:unit-obstruction}
is nonzero. Hence $a_{152}^{\mu_{h}}\ne0$, and the injectivity of
$\mathrm{d}\colon F_{4}\to F_{3}$ gives $[e_{1},e_{5},e_{2}]_{\mu_{h}}\ne0$.

For the completed ring, the same proof uses $B=\Bbbk\llbracket y,z\rrbracket$.
This is again a local domain, and every step is unchanged. 
\end{proof}
\begin{rem}
\label{rem:no-hidden} Exactness is used before reduction modulo $\langle x,w\rangle$,
to normalize the homotopy and to identify scalar associators. The
reduced complex is not assumed to remain a resolution over $B$. After
reduction, the proof uses only the displayed matrices, linear maps,
and Nakayama's lemma. In particular, no exactness assertion is being
transported through a non-flat quotient.

Likewise, no step assumes $e_{i}\star e_{i}=0$ or $e_{j}\star e_{i}=-e_{i}\star e_{j}$.
Alternation occurs only in the fixed leading products $p_{ij}$ of
the reference multiplication. The perturbations $\overline{e_{i}\star e_{j}}-p_{ij}$
and the left and right linear maps remain unrestricted. 
\end{rem}

\subsection{A family with higher-order reference obstructions}

\label{sec:family} The argument is not specific to the coefficient
$yz$ in \eqref{eq:old-associator}. Let 
\[
0\ne c(y,z)\in\langle y,z\rangle\Bbbk[y,z],\qquad I_{c}=\langle x^{2},w^{2},zw,xy,yz\,c(y,z)\rangle.
\]
The polynomial $c$ need not be homogeneous or a monomial. Put $A=R_{\mathfrak{m}}$,
with $R=\Bbbk[x,y,z,w]$. We state minimality over $A$; the construction
also provides a polynomial free resolution over $R$.
\begin{thm}
\label{thm:family} The minimal free resolution of $A/I_{c}A$ over
$A$ admits no associative multiplication. More precisely, associativity
on the $32$ core triples forces non-associativity at $(e_{1},e_{5},e_{2})$.
The same assertions hold for the explicit polynomial free resolution
over $R$ constructed below and after completion of $A$ at its maximal
ideal. 
\end{thm}

\begin{proof}
Let $J=\langle x^{2},w^{2},zw,xy\rangle$ and $f=yzc$. We first verify
$J:f=\langle x,w\rangle$. Certainly $xf$ and $wf$ belong to $J$.
Conversely, if $gf\in J$, then reduction modulo $\langle x,w\rangle$
gives $\bar{g}\,yzc=0$ in the domain $\Bbbk[y,z]$. Since $c\ne0$,
one has $\bar{g}=0$, so $g\in\langle x,w\rangle$.

The mapping cone of $R/\langle x,w\rangle\xrightarrow{\ f\ }R/J$
has the same basis as \eqref{eq:modules}. The core differentials
are unchanged, and the remaining differentials are 
\begin{equation}
\mathrm{d}e_{5}=yzc,\qquad\mathrm{d}e_{35}=we_{5}-yce_{3},\qquad\mathrm{d}e_{45}=xe_{5}-zce_{4},\qquad\mathrm{d}e_{345}=we_{45}-xe_{35}+ce_{34}.\label{eq:family-d}
\end{equation}
To see the comparison map explicitly, let $a,b$ be the Koszul generators
with $\mathrm{d}a=x$ and $\mathrm{d}b=w$. A lift to the Taylor complex
on $J$ is 
\[
1\longmapsto yzc,\qquad a\longmapsto zce_{4},\qquad b\longmapsto yce_{3},\qquad ab\longmapsto-ce_{34}.
\]
The usual cone differential, after choosing signs for its shifted
basis, gives \eqref{eq:family-d}. Since $c\in\langle y,z\rangle$,
every matrix entry is in $\mathfrak{m}$, and the localized cone is
minimal. If $c$ is homogeneous, the polynomial cone is the minimal
graded resolution over $R$ as well.

To construct a commutative reference multiplication, first work over
$\Bbbk[x,y,z,w,t]$ with ideal $\langle x^{2},w^{2},zw,xy,yzt\rangle$.
In \eqref{eq:pi}, replace each occurrence of $yz$ by $t$, each
occurrence of $y^{2}z$ by $yt$, and each occurrence of $yz^{2}$
by $zt$. These formulas define a Taylor chain retraction and hence
a commutative multiplication. Substitute $t=c(y,z)$. Chain-map and
product identities are polynomial identities and therefore survive
this substitution; exactness of the specialized complex was established
separately by the cone construction.

The reference product is still associative on the first four generators,
its top table is still \eqref{eq:top-table}, and 
\[
e_{1}e_{5}=zce_{14}+xe_{45},\qquad e_{2}e_{5}=yce_{23}+we_{35},\qquad[e_{1},e_{5},e_{2}]_{\mu}=-c\,\mathrm{d}e_{1234}.
\]
Modulo $\langle x,w\rangle$, the last entry of \eqref{eq:reduced-d}
becomes $c\bar{e}_{34}$. Since $c\in\mathfrak{n}$, the image of
the reduced differential is still contained in $\mathfrak{n}G$. The
Nakayama argument is unchanged. Assuming associativity on $\mathcal{C}$,
the final scalar associator is 
\[
\bar{a}_{152}^{\mu_{h}}=-c\bigl(1-z\bar{\ell}'_{2}(\bar{e}_{14})-y\bar{\ell}_{1}(\bar{e}_{23})\bigr)\ne0.
\]
The factor in parentheses is a unit and $c$ is nonzero in $\Bbbk[y,z]_{\langle y,z\rangle}$,
or in $\Bbbk\llbracket y,z\rrbracket$ for the completed case. Localization
and completion preserve exactness by flatness. Every polynomial multiplication
localizes, so the same relative assertion holds on the displayed polynomial
complex. This proves the theorem. 
\end{proof}
\begin{cor}
\label{cor:powers} For integers $a,b\ge1$ with $a+b\ge3$, put 
\[
I_{a,b}=\langle x^{2},w^{2},zw,xy,y^{a}z^{b}\rangle.
\]
The quotient $R/I_{a,b}$ has a minimal free resolution over $R$
with Betti numbers $(1,5,8,5,1)$ which admits no associative multiplication.
For the reference product constructed above, 
\[
[e_{1},e_{5},e_{2}]_{\mu}=-y^{a-1}z^{b-1}\mathrm{d}e_{1234}.
\]
\end{cor}

\begin{proof}
Apply Theorem~\ref{thm:family} with $c=y^{a-1}z^{b-1}$. 
\end{proof}
For example, when $a=b=r+1$ with $r\ge1$, the displayed associator
has coefficients of order $2r+1$ in $\mathfrak{m}$. Thus the reference
obstruction can lie arbitrarily deep in the maximal ideal while nonexistence
of an associative multiplication remains valid. The case $a=b=1$
is deliberately excluded from the stated complex: then the coefficient
of $e_{34}$ in $\mathrm{d}e_{345}$ is a unit and this cone is not
minimal. Cancelling that pair gives the smaller cellular resolution
appearing in Avramov's example, studied in Subsection~\ref{sec:avramov-uniform}.

\subsection{Vanishing of the regular-sequence Avramov obstructions}

\label{subsec:katthan-avramov}

Katthän noted that he could not find a regular sequence in $I$ for
which Avramov's obstruction proves nonexistence of a minimal DG algebra
resolution \cite[Remark~5.4]{Kat19}. We show that every such obstruction
vanishes.

Recall that, for a regular sequence generating $J\subseteq I$, the
obstruction associated to $\pi\colon R\to S=R/J$ and $M=R/I$ is
\[
o^{\pi}(M)=\ker\!\left(\frac{\operatorname{Tor}^{R}(M,\Bbbk)}{\operatorname{Tor}_{+}^{R}(S,\Bbbk)\operatorname{Tor}^{R}(M,\Bbbk)}\longrightarrow\operatorname{Tor}^{S}(M,\Bbbk)\right).
\]
Here the action in the denominator is the natural Cartan--Eilenberg
product; see \cite[Section~1 and Remark~1.8.1]{Avr81}.
\begin{thm}[Vanishing for all regular sequences]
\label{thm:katthan-avramov} Let $\Bbbk$ be any field, let $R=\Bbbk[x,y,z,w]$,
and put $I=\langle x^{2},w^{2},zw,xy,y^{2}z^{2}\rangle$ and $\mathfrak{m}=\langle x,y,z,w\rangle$.
Let $Q$ be $R$, $R_{\mathfrak{m}}$, or $\Bbbk\llbracket x,y,z,w\rrbracket$,
and set $M=Q/IQ$. Let $f_{1},\ldots,f_{r}$ be a $Q$-regular sequence
in $IQ$, put $J=\langle f_{1},\ldots,f_{r}\rangle$ and $S=Q/J$.
Then the canonical map 
\[
\frac{\operatorname{Tor}^{Q}(M,\Bbbk)}{\operatorname{Tor}_{+}^{Q}(S,\Bbbk)\operatorname{Tor}^{Q}(M,\Bbbk)}\longrightarrow\operatorname{Tor}^{S}(M,\Bbbk)
\]
is injective. Thus $o^{\pi}(M)=0$ in every homological degree. 
\end{thm}

\begin{proof}
For $Q=R$, localization at $\mathfrak{m}$ preserves all the Tor
groups and natural maps in the statement. We may therefore assume
that $Q$ is local, retaining the notation $Q,S,M$ after localization.
Write $K_{M}=M\otimes_{Q}K_{Q}(x,y,z,w)$, with exterior generators
$\varepsilon_{1},\ldots,\varepsilon_{4}$ and $\mathrm{d}\varepsilon_{1}=x$,
$\mathrm{d}\varepsilon_{2}=y$, $\mathrm{d}\varepsilon_{3}=z$, $\mathrm{d}\varepsilon_{4}=w$.
Write $\varepsilon_{i_{1}\cdots i_{s}}=\varepsilon_{i_{1}}\cdots\varepsilon_{i_{s}}$
for $i_{1}<\cdots<i_{s}$. We first construct a subalgebra $H\subseteq K_{M}$
of cycles representing all of its homology. Set $H_{0}=\Bbbk$ and
use these bases: 
\begin{equation}
\begin{array}{c|l}
i & \text{basis of }H_{i}\\
\hline 1 & x\varepsilon_{1},\ w\varepsilon_{4},\ w\varepsilon_{3},\ x\varepsilon_{2},\ y^{2}z\varepsilon_{3}\\[2pt]
2 & xw\varepsilon_{14},\ xw\varepsilon_{13},\ x\varepsilon_{12},\ w\varepsilon_{34},\ xw\varepsilon_{24},\ xw\varepsilon_{23},\ y^{2}z\varepsilon_{34},\ yz^{2}\varepsilon_{12}\\[2pt]
3 & xw\varepsilon_{134},\ xw\varepsilon_{124},\ xw\varepsilon_{123},\ xw\varepsilon_{234},\ y^{2}z\varepsilon_{134}\\[2pt]
4 & xw\varepsilon_{1234}.
\end{array}\label{eq:katthan-koszul-cycles}
\end{equation}
Each displayed element is a cycle, and their span is closed under
multiplication, using $x^{2}=w^{2}=zw=xy=y^{2}z^{2}=0$ in $M$. To
check that their classes are nonzero, first work over $R$ and use
the standard monomial--exterior basis in each displayed multidegree.
Every column of the incoming Koszul differential has either no nonzero
entries or two entries $1$ and $-1$. Hence the functional that sums
coefficients kills every boundary and takes value $1$ on the displayed
monomial cycle. The multidegrees are distinct within each homological
degree. The dimensions $(1,5,8,5,1)$ agree with the ranks in \eqref{eq:modules},
so these classes form a basis of $\operatorname{Tor}^{R}(R/I,\Bbbk)$.
Localization and completion are flat and preserve this Koszul homology,
which is annihilated by $\mathfrak{m}$. Thus in all three cases the
inclusion 
\[
H\hookrightarrow K_{M}
\]
is a quasi-isomorphism of DG $\Bbbk$-algebras, with zero differential
on $H$.

We use the canonical divided powers on positive even elements of the
Koszul algebra. Explicitly, if $b=\sum_{i<j}b_{ij}\varepsilon_{ij}\in(K_{M})_{2}$,
then 
\begin{equation}
\gamma_{2}(b)=(b_{12}b_{34}-b_{13}b_{24}+b_{14}b_{23})\varepsilon_{1234},\qquad\gamma_{n}(b)=0\quad(n\ge3).\label{eq:katthan-koszul-dp}
\end{equation}
These formulas satisfy $\mathrm{d}\gamma_{n}(b)=(\mathrm{d}b)\gamma_{n-1}(b)$;
in characteristic zero they agree with $b^{n}/n!$, but require no
division in arbitrary characteristic. Each individual even basis element
in \eqref{eq:katthan-koszul-cycles} has zero divided powers of order
at least two. The addition formula for divided powers and multiplicative
closure therefore show that $H$ is closed under these operations.
Its inclusion in $K_{M}$ preserves divided powers.

Write $\boldsymbol{x}=(x_{1},x_{2},x_{3},x_{4})=(x,y,z,w)$. Since
$IQ\subseteq(\mathfrak{m}Q)^{2}$, write $f_{i}=\sum_{j=1}^{4}a_{ij}x_{j}$
with $a_{ij}\in\mathfrak{m}Q$. In $K_{S}=S\otimes_{Q}K_{Q}(\boldsymbol{x})$,
the cycles $z_{i}=\sum_{j}\bar{a}_{ij}\varepsilon_{j}$ represent
the exterior generators of $\operatorname{Tor}^{Q}(S,\Bbbk)$: under
$H_{1}(K_{S})\cong J/\mathfrak{m}J$, the class $[z_{i}]$ corresponds
to $f_{i}\bmod\mathfrak{m}J$, since its lift to $K_{Q}(\boldsymbol{x})$
has differential $f_{i}$. The Tate resolution of $\Bbbk$ over the
complete intersection $S$ is 
\[
L=K_{S}\langle\delta_{1},\ldots,\delta_{r}\rangle,\qquad|\delta_{i}|=2,\qquad\mathrm{d}\delta_{i}=z_{i},
\]
where angle brackets denote divided-power variables, with $\delta_{i}^{(n)}=\gamma_{n}(\delta_{i})$;
see \cite[Section~6.1]{Avr98}. The inclusion of the Koszul complex
into $L$, after tensoring with $M$, induces the canonical change-of-rings
map.

Now work in $M\otimes_{S}L=K_{M}\langle\delta_{1},\ldots,\delta_{r}\rangle$.
Choose $h_{i}\in H_{1}$ representing $[z_{i}]$, so that $z_{i}-h_{i}=\mathrm{d}b_{i}$
for some $b_{i}\in(K_{M})_{2}$. The divided-power change of variables
$\delta_{i}'=\delta_{i}-b_{i}$ gives $\mathrm{d}\delta_{i}'=h_{i}$.
More explicitly, 
\[
(\delta_{i}')^{(n)}=\sum_{q=0}^{n}(-1)^{q}\gamma_{q}(b_{i})\delta_{i}^{(n-q)}.
\]
The divided-power identities make this an invertible DG algebra change
of variables; \eqref{eq:katthan-koszul-dp} makes it valid in every
characteristic. It fixes $K_{M}$, so the canonical change-of-rings
map is unchanged. Notice that this substitution is made after tensoring
with $M$; the elements $b_{i}$ need not lift to primitives over
$S$.

Consider the resulting morphism 
\begin{equation}
C=H\langle\delta_{1}',\ldots,\delta_{r}'\rangle\longrightarrow K_{M}\langle\delta_{1}',\ldots,\delta_{r}'\rangle,\qquad\mathrm{d}\delta_{i}'=h_{i}.\label{eq:katthan-tate-comparison}
\end{equation}
Filter both sides increasingly by total divided-power degree: $\mathcal{F}_{p}$
is spanned by terms $a(\delta_{1}')^{(n_{1})}\cdots(\delta_{r}')^{(n_{r})}$
with $\sum_{i}n_{i}\le p$. The coefficient differential preserves
this degree, and differentiating a divided-power variable lowers it
by one. On the associated graded, \eqref{eq:katthan-tate-comparison}
is therefore a direct sum of shifts of $H\hookrightarrow K_{M}$,
hence a quasi-isomorphism. Every divided-power variable has homological
degree two, so in homological degree $N$ the filtration stops at
$p=\lfloor N/2\rfloor$. The filtrations are exhaustive and bounded
in each degree, so the comparison of their spectral sequences proves
that \eqref{eq:katthan-tate-comparison} is a quasi-isomorphism.

Since $\mathrm{d}H=0$, the differential of $C$ lowers divided-power
degree exactly one. The degree-zero divided-power component of a boundary
consequently belongs to $\sum_{i}h_{i}H$. Conversely, each $h_{i}a$
is a boundary, by differentiating $\delta_{i}'a$. It follows that
\[
\ker\bigl(H\longrightarrow H(C)\bigr)=\sum_{i}h_{i}H.
\]
Under $H\cong\operatorname{Tor}^{Q}(M,\Bbbk)$, the elements $h_{i}$
are the images of the degree-one generators of $\operatorname{Tor}^{Q}(S,\Bbbk)$.
Since the latter algebra is exterior on these generators, $\sum_{i}h_{i}H$
is exactly $\operatorname{Tor}_{+}^{Q}(S,\Bbbk)\operatorname{Tor}^{Q}(M,\Bbbk)$.
The preceding kernel equality proves the asserted injectivity. 
\end{proof}
Together with Theorem~\ref{thm:main}, this gives a monomial example
with no associative multiplication on its minimal resolution over
any field, although every Avramov obstruction arising from a regular
sequence in its defining ideal vanishes.

\section{Obstructions to multiplicative Taylor comparisons}

\label{sec:tcomp}

Throughout this section, let $\Bbbk$ be a field, let $R=\Bbbk[x_{1},\ldots,x_{n}]$,
let $\mathfrak{m}=\langle x_{1},\ldots,x_{n}\rangle$, let $I\subsetneq R$
be a nonzero monomial ideal, and let $F$ be the minimal multigraded
free resolution of $R/I$ over $R$, with $F_{0}=R$. The resolution
carries its fine multigrading, but competing multiplications and comparison
maps need not preserve it.

Let $J\subseteq I$ be generated by a selection of minimal generators
of $I$, and equip its Taylor resolution $T(J)$ with the usual associative
multiplication $\mu_{T}$ \cite[Section 1.2]{Kat19}. A unital chain
comparison $\varphi\colon T(J)\to F$ is an $R$-linear chain map
with $\varphi_{0}=\operatorname{id}_{R}$, lifting $R/J\twoheadrightarrow R/I$.
It is multiplicative for a multiplication $\mu_{h}$ on $F$ if

\[
\varphi(ab)=\mu_{h}(\varphi(a),\varphi(b))\qquad(a,b\in T(J)).
\]
We ask how large a selection can admit such a comparison when the
multiplication on $F$ is allowed to vary.

Uniform non-associativity gives an obstruction that survives changes
of both the multiplication and the degree-one lifts defining the comparison.
For the edge ideal of the eight-cycle, the largest possible selection
has four generators. Its big height is five, and an injective comparison
on five selected generators exists, but every five-generator selection
fails to admit a multiplicative comparison. We also obtain a family
of cycle edge ideals with this failure at the big-height size.

For squarefree ideals, Katthän realizes a multigraded DG algebra resolution
as a quotient of the full Taylor algebra \cite[Theorem 3.6]{Kat19}.
VandeBogert extends this to generalized Taylor resolutions \cite[Section 4]{VdB23},
and Ferraro and Alvarez give a divided-power version \cite[Theorem 5.4]{FA25}.
Those results assume an associative multigraded target. Our question
concerns Taylor algebras on selected generators and permits a nonassociative
target.

This question is motivated by multiplicative approaches to rank bounds.
Buchsbaum and Eisenbud use an associative commutative multiplication
on a resolution to obtain an injective comparison from the Koszul
algebra on a regular sequence \cite[Proposition 1.4]{BE77}. Here
the source is a Taylor algebra on selected minimal generators, and
the target multiplication is allowed to be nonassociative.

\subsection{Injective Taylor comparisons and rank bounds}

We first record an injectivity consequence of the associated-prime
characterization by dominant sets in \cite[Theorem 2.1]{BEK+25}.
That result supplies the selected generators and the pure-power ideal
used below; we include the construction to make the comparison-map
argument explicit. The proposition distinguishes failure of multiplicativity
from a failure of the underlying numerical bound.
\begin{prop}[Injective Taylor comparisons]
\label{prop:tcomp-injective} Let $I\subsetneq R$ be a nonzero monomial
ideal, and let 
\[
m=\operatorname{bight}I:=\max\{\mathrm{ht}\,\mathfrak{p}\mid\mathfrak{p}\in\operatorname{Ass}_{R}(R/I)\}.
\]
There are $m$ distinct minimal generators $u_{1},\ldots,u_{m}$ of
$I$ such that $R/J$, where $J=\langle u_{1},\ldots,u_{m}\rangle$,
has a minimal Taylor resolution $T$ over $R$, and every unital multigraded
comparison $\varphi\colon T\to F$ lifting $R/J\twoheadrightarrow R/I$
is injective. In particular, 
\[
\beta_{i}^{R}(R/I)\geq\binom{m}{i}.
\]
\end{prop}

\begin{proof}
Choose $\mathfrak{p}\in\operatorname{Ass}_{R}(R/I)$ with $\mathrm{ht}\,\mathfrak{p}=m$.
Since $I$ is monomial, we can write $\mathfrak{p}=\langle z_{1},\ldots,z_{m}\rangle$,
where $z_{j}=x_{i_{j}}$ for distinct indices $i_{1},\ldots,i_{m}\in\{1,\ldots,n\}$.
Thus $z_{1},\ldots,z_{m}$ is an ordered selection of the original
variables, and $m\leq n$. Choose a monomial $v\notin I$ such that
$I:v=\mathfrak{p}$, and put $a_{j}=\deg_{z_{j}}(v)+1$. Since $z_{j}v\in I$,
choose a minimal generator $u_{j}$ of $I$ dividing $z_{j}v$. Since
$u_{j}\nmid v$, we have 
\begin{equation}
\deg_{z_{j}}(u_{j})=a_{j},\qquad\deg_{z_{j'}}(u_{j})<a_{j'}\quad(j'\ne j).\label{eq:tcomp-private}
\end{equation}
The selected generators are distinct. For each $j$, the $z_{j}$-exponent
of $u_{j}$ is strictly larger than that of every other selected generator.
Consequently every coefficient in the positive Taylor differential
belongs to $\mathfrak{m}$. Hence the Taylor resolution of $R/J$
over $R$ is minimal. The selected generators form a dominant set
in the terminology of \cite[Definition 4.1]{Ale17}, whose Theorem
4.4 characterizes minimal Taylor resolutions by this condition.

Set $Q=\langle z_{1}^{a_{1}},\ldots,z_{m}^{a_{m}}\rangle$. We claim
that $I\subseteq Q$. For a minimal generator $u$ of $I$, the monomial
$r=u/\gcd(u,v)$ belongs to $I:v=\mathfrak{p}$. Thus some $z_{j}$
divides $r$, which means $\deg_{z_{j}}(u)\geq a_{j}$. This proves
the claim.

Let $E=K_{R}(z_{1}^{a_{1}},\ldots,z_{m}^{a_{m}})$, and choose a unital
multigraded comparison $\psi\colon F\to E$ lifting $R/I\twoheadrightarrow R/Q$.
The inclusions $J\subseteq I\subseteq Q$ give chain comparisons
\[
T\xrightarrow{\varphi}F\xrightarrow{\psi}E.
\]
We show that their composition is injective. For $A\subseteq\{1,\ldots,m\}$,
put 
\[
u_{A}=\operatorname{lcm}(u_{j}:j\in A),\qquad q_{A}=\prod_{j\in A}z_{j}^{a_{j}}.
\]
The Taylor algebra $T$ and the Koszul algebra $E$ are both supported
on the full simplex on $\{1,\ldots,m\}$, with respective face labels
$u_{A}$ and $q_{A}$. Write $t_{A}$ and $\varepsilon_{A}$ for the
Taylor and Koszul basis elements, respectively. A coefficient from
$t_{A}$ to $\varepsilon_{B}$, with $|A|=|B|$, can be nonzero only
if $q_{B}\mid u_{A}$. By \eqref{eq:tcomp-private}, this forces $B\subseteq A$
and hence $B=A$. Thus every unital multigraded chain map $\theta\colon T\to E$
is diagonal. Commuting with the differentials, starting with $\theta_{0}=1$,
gives 
\[
\theta(t_{A})=\frac{u_{A}}{q_{A}}\varepsilon_{A}.
\]
All diagonal coefficients are nonzero monomials. Taking $\theta=\psi\varphi$,
we see that $\psi\varphi$ is injective, and hence $\varphi$ is injective.
The rank inequality follows. 
\end{proof}
For Katthän's ideal $I=\langle x^{2},w^{2},zw,xy,y^{2}z^{2}\rangle$,
take $v=xw$. Then $v\notin I$ and $I:v=\langle x,w,z,y\rangle=\mathfrak{m}$,
so $m=n=4$. The construction gives
\[
\begin{gathered}(z_{1},z_{2},z_{3},z_{4})=(x,w,z,y),\qquad(a_{1},a_{2},a_{3},a_{4})=(2,2,1,1),\\
(u_{1},u_{2},u_{3},u_{4})=(x^{2},w^{2},zw,xy),\qquad Q=\langle x^{2},w^{2},z,y\rangle.
\end{gathered}
\]
The exponents are strictly largest among these four selected generators;
the omitted generator $y^{2}z^{2}$ has larger $y$- and $z$-exponents.
In degree one, the diagonal comparison sends $(t_{1},t_{2},t_{3},t_{4})$
to $(\varepsilon_{1},\varepsilon_{2},w\varepsilon_{3},x\varepsilon_{4})$.
For the reference multiplication in Subsection \ref{sec:reference},
this four-generator Taylor algebra is a subalgebra of $F$.

For comparison, put $c=\mathrm{ht}\,I$ and $\rho=\operatorname{pd}_{R}(R/I)$.
Then
\[
c\leq m=\operatorname{bight}I\leq\rho.
\]
The last inequality follows from the Auslander--Buchsbaum formula
after localization at an associated prime. Brun and Römer prove the
stronger numerical bound
\[
\beta_{i}^{R}(R/I)\geq\binom{\rho}{i}
\]
for monomial quotients over every field \cite[Corollary 4.1(i)]{BR04}.
For related binomial bounds from dominant sets, see \cite[Theorem 3.5]{Ale23}.
Proposition \ref{prop:tcomp-injective} supplies specific injective
comparisons to which we can apply the question of multiplicativity.

\subsection{Changing lifts and the induced-path obstruction}

The obstruction must be independent of the choice of comparison. The
following elementary consequence of minimality is useful for this
purpose.
\begin{lem}[Changing the degree-one lifts]
\label{lem:tcomp-lifts} Let $(R,\mathfrak{m})$ be local and let
$F$ be a minimal free resolution of $R/I$ over $R$ with $F_{0}=R$,
and let $\mu$ be a multiplication on $F$. Then
\[
[F,F,F]_{\mu}\subseteq\mathfrak{m}F.
\]
If $a_{j},a'_{j}\in F_{1}$ and $a'_{j}-a_{j}\in\mathfrak{m}F_{1}$
for $1\leq j\leq3$, then 
\begin{equation}
[a'_{1},a'_{2},a'_{3}]_{\mu}\equiv[a_{1},a_{2},a_{3}]_{\mu}\pmod{\mathfrak{m}^{2}F_{3}}.\label{eq:tcomp-lift-change}
\end{equation}
In particular, an associator which is nonzero modulo $\mathfrak{m}^{2}F_{3}$
cannot be removed by changing these degree-one lifts. 
\end{lem}

\begin{proof}
Put $\Bbbk=R/\mathfrak{m}$ and $F_{\Bbbk}=F\otimes_{R}\Bbbk$. By
minimality, the differential on $F_{\Bbbk}$ vanishes, so
\[
F_{\Bbbk}=H(F_{\Bbbk})\cong\operatorname{Tor}^{R}(R/I,\Bbbk).
\]
The product induced by $\mu$ is the usual associative product on
this Tor algebra. Indeed, the two iterated products $\mu(\mu\otimes1)$
and $\mu(1\otimes\mu)$ on $F$ lift the same triple product on $R/I$,
so they are chain homotopic. Tensoring this homotopy with $\Bbbk$
makes the two maps equal because the differentials vanish. Thus tensoring
with $\Bbbk$ kills every associator, giving $[F,F,F]_{\mu}\subseteq\mathfrak{m}F$.

For $z\in\mathfrak{m}$ and $a,b,c\in F_{1}$, trilinearity and the
first assertion give
\[
[za,b,c]_{\mu}=z[a,b,c]_{\mu}\in\mathfrak{m}^{2}F_{3}.
\]
Taking finite sums gives the same conclusion whenever the first input
lies in $\mathfrak{m}F_{1}$, and the same argument applies in the
other two positions. Since $a'_{j}-a_{j}\in\mathfrak{m}F_{1}$, replacing
one input at a time yields
\[
\begin{aligned}{}[a'_{1},a'_{2},a'_{3}]_{\mu} & \equiv[a_{1},a'_{2},a'_{3}]_{\mu}\\
 & \equiv[a_{1},a_{2},a'_{3}]_{\mu}\\
 & \equiv[a_{1},a_{2},a_{3}]_{\mu}\pmod{\mathfrak{m}^{2}F_{3}}.
\end{aligned}
\]
This proves \eqref{eq:tcomp-lift-change}.
\end{proof}
Suppose now that $T$ is a Taylor resolution over $R$ on a selection
of minimal generators of $I$, and $t_{i}\in T_{1}$ and $e_{i}\in F_{1}$
have the same differential $u_{i}$. For every unital chain comparison
$\varphi\colon T\to F$, 
\begin{equation}
\varphi(t_{i})-e_{i}\in\ker(\mathrm{d}\colon F_{1}\to F_{0})=\mathrm{d}F_{2}\subseteq\mathfrak{m}F_{1}.\label{eq:tcomp-comparison-lift}
\end{equation}
Thus Lemma~\ref{lem:tcomp-lifts} applies to the images under every
comparison, not only to a fixed multigraded comparison.
\begin{prop}
\label{prop:tcomp-path} Let $\Gamma$ be a finite simple graph, let
$I(\Gamma)$ be its edge ideal, and let $F$ be the minimal multigraded
free resolution of $R/I(\Gamma)$ over $R$. Suppose $\Gamma$ contains
an induced path 
\[
v_{1}-v_{2}-v_{3}-v_{4}-v_{5}-v_{6}.
\]
Let $J$ be generated by a selection of edges of $\Gamma$ which includes
\[
v_{1}v_{2},\qquad v_{3}v_{4},\qquad v_{5}v_{6},
\]
and give $T(J)$ its usual Taylor multiplication. Then no unital chain
comparison $T(J)\to F$ lifting $R/J\twoheadrightarrow R/I(\Gamma)$
is multiplicative for any multiplication on $F$. Neither the comparison
nor the competing multiplication is required to respect an internal
grading. 
\end{prop}

\begin{proof}
The path calculation in Subsection~\ref{subsec:path} satisfies the
hypotheses of Theorem~\ref{theoremnonassociativeforall}. Since the
path is induced, its edge ideal is an elimination ideal of $I(\Gamma)$.
Proposition~\ref{prop:elimination-persistence} transfers those hypotheses
to $F$ at the degree-one elements $a_{1},a_{2},a_{3}$ corresponding
to the three displayed alternating edges.

The proof of the uniform criterion gives more than nonvanishing: for
every multiplication $\nu$, 
\[
[a_{1},a_{2},a_{3}]_{\nu}\notin\mathfrak{m}^{2}F_{3}.
\]
Indeed, if this associator belonged to $\mathfrak{m}^{2}F_{3}$, the
variation congruence \eqref{eq:uniform-congruence} would place the
reference associator in the submodule excluded by the criterion. The
same reasoning applies after the elimination step. One may localize
at the homogeneous maximal ideal throughout.

Let $t_{1},t_{2},t_{3}\in T(J)_{1}$ denote the corresponding Taylor
generators. Equations \eqref{eq:tcomp-comparison-lift} and \eqref{eq:tcomp-lift-change}
imply 
\[
[\varphi(t_{1}),\varphi(t_{2}),\varphi(t_{3})]_{\nu}\ne0.
\]
But multiplicativity and associativity of the Taylor algebra would
give 
\[
[\varphi(t_{1}),\varphi(t_{2}),\varphi(t_{3})]_{\nu}=\varphi([t_{1},t_{2},t_{3}]_{\mu_{T}})=0,
\]
a contradiction. 
\end{proof}
The Taylor resolution in Proposition~\ref{prop:tcomp-path} need
not be minimal. The obstruction rules out multiplicativity itself,
before any question about injectivity is considered.

\subsection{Cycle edge ideals}

The cycle edge ideal considered here is different from the ideal of
the eight-vertex Bier sphere in Theorem~\ref{thm:gorenstein}. Its
role is to obstruct every sufficiently large Taylor selection.
\begin{thm}
\label{thm:tcomp-c8} Let $R=\Bbbk[x_{1},\ldots,x_{8}]$ and 
\[
I=\langle u_{1},\ldots,u_{8}\rangle,\qquad u_{i}=x_{i}x_{i+1},
\]
where indices are read modulo $8$. Let $F$ be the minimal free resolution
of $R/I$ over $R$. Then $\mathrm{ht}\,I=4$ and $\operatorname{bight}I=5$.
For every five-element subset $S\subseteq\{1,\ldots,8\}$, put $J_{S}=\langle u_{i}:i\in S\rangle$
and let $T_{S}$ be its Taylor algebra. No unital comparison $T_{S}\to F$
is multiplicative for any multiplication on $F$.

Nevertheless, a five-generator selection has a minimal Taylor resolution
over $R$ admitting a degreewise injective comparison into $F$, and
a four-generator selection admits a multiplicative embedding for a
suitable multigraded multiplication on $F$. Thus the largest number
of selected minimal generators whose Taylor algebra can map multiplicatively
into $F$ is exactly four. These assertions hold over every field. 
\end{thm}

\begin{proof}
The four disjoint edges $12,34,56,78$ show that every vertex cover
has at least four vertices, while the even vertices form a cover of
size four. Thus $\mathrm{ht}\,I=4$.

The complements of maximal independent sets are the minimal vertex
covers. Every maximal independent set of $C_{8}$ is dominating, and
a vertex together with its neighbors contains at most three vertices.
Hence every maximal independent set has at least three vertices. The
set $\{2,5,7\}$ is a maximal independent set of size three. The largest
minimal vertex cover therefore has size five, giving $\operatorname{bight}I=5$.

Now let $S$ have five elements. At least three elements of $S$ have
the same parity. Every three-element subset of either $\{1,3,5,7\}$
or $\{2,4,6,8\}$ has the form 
\[
\{i,i+2,i+4\}\pmod 8.
\]
The corresponding generators are the alternating edges of the induced
path on the six consecutive vertices $i,i+1,\ldots,i+5$. Proposition~\ref{prop:tcomp-path}
proves the nonexistence assertion.

To see that this is not a failure of injectivity, take 
\[
\begin{aligned}J_{5} & =\langle u_{1},u_{2},u_{4},u_{5},u_{7}\rangle\\
 & =\langle x_{1}x_{2},x_{2}x_{3},x_{4}x_{5},x_{5}x_{6},x_{7}x_{8}\rangle.
\end{aligned}
\]
Its generators have private variables $x_{1},x_{3},x_{4},x_{6},x_{8}$,
respectively, so its Taylor resolution is minimal. Moreover, 
\[
I:\langle x_{2}x_{5}x_{7}\rangle=\langle x_{1},x_{3},x_{4},x_{6},x_{8}\rangle.
\]
These are precisely the witness and selection used in Proposition~\ref{prop:tcomp-injective}.
Thus the multigraded comparison $T(J_{5})\to F$ is injective, although
no comparison from this Taylor algebra can be multiplicative.

Finally, let 
\[
\begin{aligned}J_{4} & =\langle u_{1},u_{2},u_{5},u_{6}\rangle\\
 & =\langle x_{1}x_{2},x_{2}x_{3},x_{5}x_{6},x_{6}x_{7}\rangle.
\end{aligned}
\]
Put $R_{0}=\Bbbk[x_{1},x_{2},x_{3},x_{5},x_{6},x_{7}]$. Then $J_{4}=(I\cap R_{0})R$.
The Taylor resolution of $R_{0}/(I\cap R_{0})$ over $R_{0}$ is minimal,
with private variables $x_{1},x_{3},x_{5},x_{7}$, and its scalar
extension is the Taylor resolution $T(J_{4})$ of $R/J_{4}$ over
$R$. Lemma~\ref{lem:specialization-retraction} identifies this
scalar extension with a basis subcomplex of $F$ and extends its Taylor
multiplication to a multigraded multiplication on $F$. The resulting
inclusion is multiplicative and injective. The same positive conclusion
also follows from \cite[Proposition 3.1]{Kat19}: the four selected
generators span a full Scarf simplex, since they have private variables
and every omitted edge uses $x_{4}$ or $x_{8}$. A multiplicative
comparison on more than five selected generators would restrict to
one on five, which has already been excluded. This proves the exact
maximum. 
\end{proof}
\begin{figure}[htbp]
\centering %
\begin{tabular}{cc}
\begin{tikzpicture}[scale=0.9]
\foreach \i/\a in {1/90,2/45,3/0,4/-45,5/-90,6/-135,7/180,8/135}
  \node[circle,fill=black,inner sep=1.5pt,label=\a:$\i$] (v\i) at (\a:1.6) {};
\foreach \i/\j in {1/2,2/3,3/4,4/5,5/6,6/7,7/8,8/1}
  \draw[gray,densely dashed] (v\i)--(v\j);
\foreach \i/\j in {1/2,2/3,4/5,5/6,7/8}
  \draw[black,thick] (v\i)--(v\j);
\foreach \i/\j in {1/2,5/6,7/8}
  \draw[red!75!black,line width=1.6pt] (v\i)--(v\j);
\node at (0,-2.70) {$J_5=\langle u_1,u_2,u_4,u_5,u_7\rangle$};
\end{tikzpicture}\qquad{} & \begin{tikzpicture}[scale=0.9]
\foreach \i/\a in {1/90,2/45,3/0,4/-45,5/-90,6/-135,7/180,8/135}
  \node[circle,fill=black,inner sep=1.5pt,label=\a:$\i$] (v\i) at (\a:1.6) {};
\foreach \i/\j in {1/2,2/3,3/4,4/5,5/6,6/7,7/8,8/1}
  \draw[gray,densely dashed] (v\i)--(v\j);
\foreach \i/\j in {1/2,2/3,5/6,6/7}
  \draw[blue!65!black,line width=1.1pt] (v\i)--(v\j);
\node at (0,-2.70) {$J_4=\langle u_1,u_2,u_5,u_6\rangle$};
\end{tikzpicture}\tabularnewline
\end{tabular}\caption{Solid edges are the selected generators. In $J_{5}$, the three heavier
edges $u_{5},u_{7},u_{1}$ alternate along the induced path $5-6-7-8-1-2$
and obstruct multiplicativity. Deleting vertices $4,8$ leaves exactly
$J_{4}$, whose Taylor algebra embeds multiplicatively.}
\label{fig:tcomp-selections} 
\end{figure}
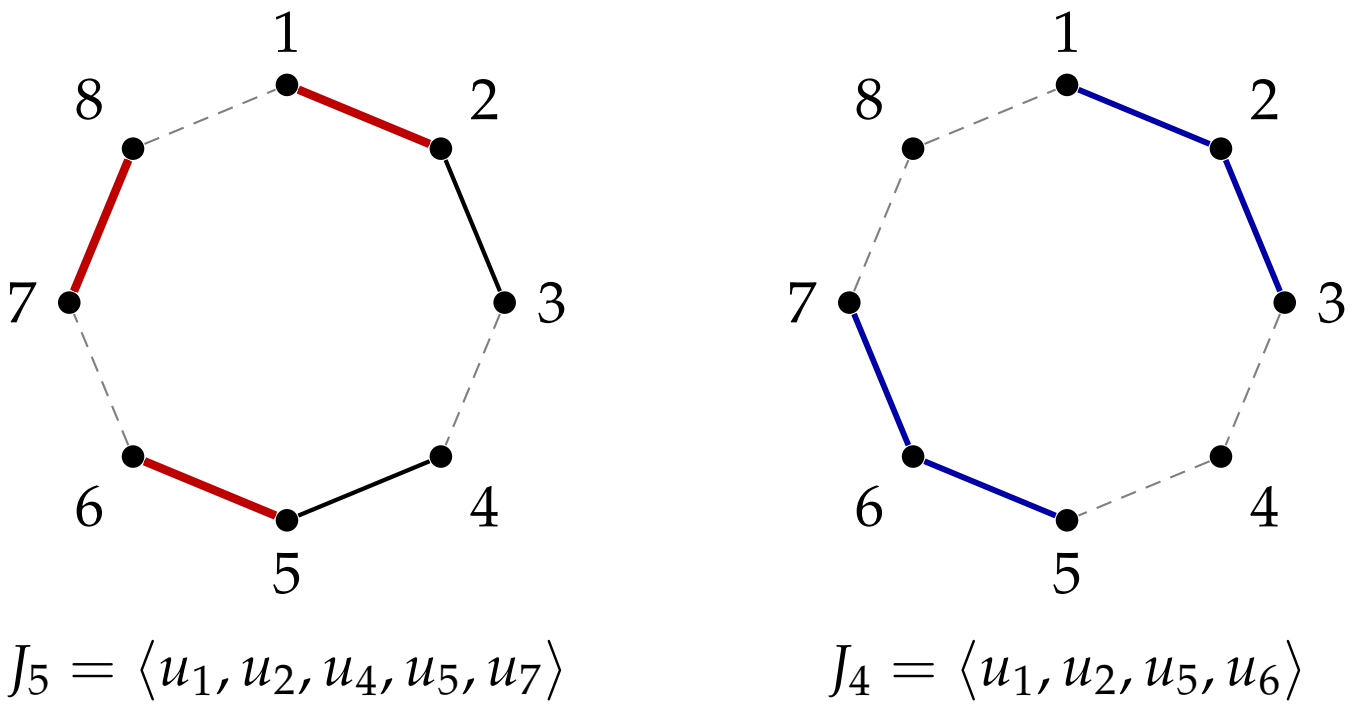

\begin{rem}
The nonexistence of a commutative DG algebra structure on the minimal
resolution of the eight-cycle is already a consequence of the classification
of DG cycles in \cite[Theorem~6.2]{GMP26}. Theorem~\ref{thm:tcomp-c8}
concerns a different question. The target is allowed to remain nonassociative,
and the multiplication, selection, and comparison may all vary. Even
with this freedom, multiplicativity cannot be arranged on five selected
generators. This contrasts with the four-generator Taylor subalgebra
in Katthän's example. 
\end{rem}

\begin{cor}
\label{cor:tcomp-cycle-family} For $r\geq1$, put $R=\Bbbk[x_{1},\ldots,x_{6r+2}]$,
let $I$ be the edge ideal of $C_{6r+2}$ and let $F$ be the minimal
resolution of $R/I$ over $R$. Then 
\[
\mathrm{ht}\,I=3r+1,\qquad\operatorname{bight}I=4r+1.
\]
For any selection of $4r+1$ edge generators, let $J$ be the ideal
they generate. No unital chain comparison $T(J)\to F$ lifting $R/J\twoheadrightarrow R/I$
is multiplicative for any multiplication on $F$. Neither the comparison
nor the multiplication is required to respect the internal grading. 
\end{cor}

\begin{proof}
A perfect matching and a parity vertex cover give height $3r+1$.
Every maximal independent set has at least $\lceil(6r+2)/3\rceil=2r+1$
vertices, and the set $\{1,4,7,\ldots,6r+1\}$ attains this bound.
Thus big height is $4r+1$.

The odd and even edge indices each form a cyclically ordered set of
size $3r+1$ when successive indices differ by two. If a selection
avoids all triples $\{i,i+2,i+4\}$, its indicator sequence on either
parity class has no three consecutive ones. Summing the $3r+1$ cyclic
inequalities 
\[
a_{j}+a_{j+1}+a_{j+2}\leq2
\]
gives $3\sum_{j}a_{j}\leq2(3r+1)$, hence $\sum_{j}a_{j}\leq2r$.
Such a selection has at most $4r$ edges. Every selection of $4r+1$
edges therefore contains the alternating edges of an induced six-vertex
path, and Proposition~\ref{prop:tcomp-path} applies. 
\end{proof}

\section*{Declaration of generative AI and AI-assisted technologies}

During the preparation of this work, the author used OpenAI Codex
to assist with exposition, literature searches, exploration and checking
of proof arguments, and the generation of symbolic-verification scripts.
The author reviewed and revised the resulting material and takes responsibility
for the mathematical claims, computations, and references.

\end{document}